\documentclass[12pt]{article}

\usepackage{amssymb,amsmath}

\newtheorem{LL}{Lemma}
\newtheorem{TT}{Theorem}
\newtheorem{CC}{Corollary}
\newtheorem{PP}{Proposition}

\newtheorem{RR}{Remark}

\newcommand{\Ad}  
{\mathbb{A}}

\newcommand{\adele} 
{\textrm{ad\`{e}le }}

\newcommand{\Adele} 
{\textrm{Ad\`{e}le }}

\newcommand{\adeles} 
{\textrm{ad\`{e}les }}

\newcommand{\alpharho} 
{\tilde{\rho}}

\newcommand{\arcsec} 
{\mathrm{arcsec}}

\newcommand{\bs} 
{\backslash}

\newcommand{\C}  
{{\cal C}}

\newcommand{\cct}  
{\mathsf{H}}

\newcommand{\cfld}  
{\mathbb{C}}

\newcommand{\CFT}  
{\textrm{class field tower}}

\newcommand{\cn}  
{{\cal CN}}

\newcommand{\code}  
{{\cal C}}

\newcommand{\Coker}  
{\mathrm{Coker}}

\newcommand{\coxeter}  
{\mathrm{coxeter}}

\newcommand{\cv}  
{{\cal V}}

\newcommand{\dfnt}  
{\stackrel{\mathrm{def}}{=}}

\newcommand{\dH}  
{d_\mathrm{H}}

\newcommand{\diag}  
{\mathrm{diag}}

\newcommand{\emb}
{\hookrightarrow}

\newcommand{\eof}  
{\hspace*{\fill}~\(\Box\)\par\endtrivlist\unskip}

\newcommand{\even}  
{\mathrm{even}}

\newcommand{\FF}  
{\mathbb{F}}

\newcommand{\Fq}  
{\mathbb{F}_q}

\newcommand{\fro}  
{\mathrm{euc}}

\newcommand{\Gal}  
{\mathrm{Gal}}

\newcommand{\idele} 
{\textrm{id\`{e}le }}

\newcommand{\Idele} 
{\textrm{Id\`{e}le }}

\newcommand{\ideles} 
{\textrm{id\`{e}les }}

\newcommand{\idx}  
{{\cal I}}

\newcommand{\iffi}  
{\textrm{if and only if }}

\newcommand{\ifld}  
{\mathbb{Z}}

\newcommand{\imag}  
{\mathrm{Im}}

\newcommand{\Image}  
{\mathrm{Im}}

\newcommand{\iso}  
{\cong}

\newcommand{\IU}  
{\boldsymbol{i}}

\newcommand{\Ker}  
{\mathrm{Ker}}

\newcommand{\lf} 
{\lfloor} 

\newcommand{\Li}  
{\mathrm{Li}}

\newcommand{\msqrt}[1]  
{\sqrt[\leftroot{-4} \uproot{2} M]{#1}}

\newcommand{\msqrts}[1]  
{\sqrt[\leftroot{-1} \uproot{2} M]{#1}}

\newcommand{\mx} 
{\mbox{ }}

\newcommand{\nequiv} 
{\equiv\!\!\!\!\!/}

\newcommand{\nfld}  
{\mathbb{N}}

\newcommand{\NKQ}  
{\mathrm{N}_{K/\mathbb{Q}}}

\newcommand{\norm}  
{\mathfrak{N}}

\newcommand{\odd}  
{\mathrm{odd}}

\newcommand{\ol}  
{\overline}

\newcommand{\Pfld}  
{\mathbb{P}}

\newcommand{\prob}  
{\Pr}

\newcommand{\Qfld}  
{\mathbb{Q}}

\newcommand{\rank}  
{\mathrm{rank}}

\newcommand{\rankin}  
{\mathrm{rankin}}

\newcommand{\real}  
{\mathrm{Re}}

\newcommand{\rf} 
{\rfloor} 

\newcommand{\rfld}  
{\mathbb{R}}

\newcommand{\Rfld}  
{\mathbb{R}}

\newcommand{\RGV}  
{R_\mathrm{GV}(\delta,q)}

\newcommand{\rt}  
{\mathsf{T}}

\newcommand{\scode}  
{\mathrm{code}}

\newcommand{\sgn}  
{\mathrm{sgn}}

\newcommand{\smsqrt}[1]  
{\sqrt[\leftroot{-2} \uproot{4} 2m]{#1}}

\newcommand{\smsqrts}[1]  
{\sqrt[\leftroot{-2} \uproot{2} 2m]{#1}}

\newcommand{\spack}  
{\mathrm{pack}}

\newcommand{\tr}  
{\mathrm{Tr}}

\newcommand{\ve}  
{\varepsilon}

\newcommand{\vol}  
{\mathrm{vol}}

\newcommand{\wt}  
{\widetilde}

\newcommand{\xpct}  
{\mathbf{E}}

\newcommand{\zfld}  
{\mathbb{Z}}

\newcommand{\Zfld}  
{\mathbb{Z}}

\begin{document}

{
\title{
\vspace{-1.3cm}
\Large
\bf  The \(q\)-ary Gilbert--Varshamov bound can be improved for all but finitely many positive integers \(q\)
}

\author{
\vspace{-0.3cm}
{\small Xue-Bin Liang}
}

\date{}

\maketitle        

\vspace{-0.6cm}

\begin{abstract}
\footnotesize
For any positive integer \(q\geq 2\) and any real number \(\delta\in(0,1)\),
let \(\alpha_q(n,\delta n)\) denote the maximum size of
a subset of \(\zfld_q^n\) with minimum Hamming distance at least \(\delta n\), where \(\zfld_q=\{0,1,\dotsc,q-1\}\) and  \(n\in\nfld\).
The asymptotic rate function is defined by
\(
R_q(\delta) = \limsup_{n\rightarrow\infty}\frac{1}{n}\log_q\alpha_q(n,\delta n).
\)
The famous \(q\)-ary asymptotic Gilbert--Varshamov bound, obtained in the 1950s, states that
\[
R_q(\delta) 
\geq 1 - \delta\log_q(q-1)-\delta\log_q\frac{1}{\delta}-(1-\delta)\log_q\frac{1}{1-\delta}
\dfnt\RGV
\]
for all positive integers \(q\geq 2\) and \(0<\delta<1-q^{-1}\). 
In the case that \(q\) is an even power of a prime with \(q\geq 49\), 
the \(q\)-ary Gilbert--Varshamov bound was firstly improved
by using algebraic geometry codes in the works of Tsfasman, Vladut, and Zink and of Ihara in the 1980s.
These algebraic geometry codes have been modified
to improve the \(q\)-ary Gilbert--Varshamov bound \(\RGV\) at a specific tangent point
\(\delta=\delta_0\in (0,1)\) of the curve \(\RGV\) for each given integer \(q\geq 46\).
However, the \(q\)-ary Gilbert--Varshamov bound \(\RGV\)
at \(\delta=1/2\), i.e., \(R_\mathrm{GV}(1/2,q)\),
remains the largest known lower bound of \(R_q(1/2)\)
for {\em infinitely many} positive integers \(q\) 
which is a generic prime and which is a generic non-prime-power integer.

In this paper, we prove that the \(q\)-ary Gilbert--Varshamov bound \(\RGV\) with \(\delta\in(0,1)\)
can be improved for all but {\em finitely many} positive integers \(q\).
Precisely, for any given number \(\ve\in(0,\frac{1}{2})\), if \(q\) is large enough,
then \(R_q(\delta)  > R_\mathrm{GV}(\delta,q)\) for all \(\delta\in[\ve,1-\ve]\).
In particular, it is explicitly shown that
\(
R_q(1/2)  > R_\mathrm{GV}(1/2,q)
\)
for all integers \(q > \exp(29)\). 
Furthermore, based on the \(q\)-ary asymptotic Plotkin upper bound
\(
    R_q(\delta) \leq 1-\delta-\frac{\delta}{q-1}
\) 
for all positive integers \(q\geq 2\) and \(0<\delta<1-q^{-1}\),
we can define the growth on the rate function \(R_q(\delta)\) for \(\delta\in(0,1)\) fixed and \(q\) growing large as follows
\[
\eta(\delta)= \liminf_{q\rightarrow\infty}\frac{1}{\log q}\log[1-\delta-R_q(\delta)]^{-1}.
\]
The \(q\)-ary Plotkin upper bound implies that the growth satisfies \(\eta(\delta)\leq 1\),
while the \(q\)-ary Gilbert--Varshamov lower bound only provides the trivial bound \(\eta(\delta)\geq 0\).
In this paper, we show that the growth has actually a nontrivial lower bound as \(\eta(\delta)\geq 1/6\) for \(\delta\in(0,1)\).
These new lower bounds, which exceed the Gilbert--Varshamov bound,
are achieved by using codes from geometry of numbers introduced by Lenstra in the 1980s.

\end{abstract}

}


\section{Introduction}
\label{sec:introduction}

We consider as an alphabet any finite set of \(q\) symbols such as \(\zfld_q=\{0,1,\dotsc,q-1\}\)
for any positive integer \(q\geq 2\).
Let \(\zfld_q^n\) denote the set of all \(n\)-tuples 
\(x=(x_1,\dotsc,x_n)\) with components in \(\zfld_q\) for any positive integer \(n\geq 1\).
We define a metric on the set \(\zfld_q^n\)
by letting \(\dH(x,y)\) be the number of positions 
in which the components of the two \(n\)-tuples \(x\) and \(y\) in \(\zfld_q^n\) differ,
which is called the Hamming distance between \(x\) and \(y\).
A \(q\)-ary code with length \(n\) is a nonempty subset \(\C\) of~\(\zfld_q^n\)
and its elements are called codewords.
Let \(\nfld\) and \(\Pfld\) denote the sets of positive integers and of prime numbers,
respectively.
If the size of \(\C\), denoted by \(|\C|\), is \(M\in\nfld\) with \(M \geq 2\)  
and the minimum Hamming distance between any two distinct codewords in \(\C\) is \(d\),
then \(\C\) is called a \(q\)-ary \((n,M,d)\) code. 
The rate \(R\) of the code \(\C\subseteq\zfld_q^n\) is defined by
\(
R=\frac{1}{n}\log_q |\C|,
\)
which is a value in \([0,1]\).
For an arbitrary real number \(\delta\in (0,1)\), let \(\alpha_q(n,\delta n)\) denote the maximum size of
a \(q\)-ary code with length \(n\) which has the minimum Hamming distance at least~\(\delta n\).
We define the asymptotic rate function as follows
\[
R_q(\delta) = \limsup_{n\rightarrow\infty}\frac{1}{n}\log_q\alpha_q(n,\delta n).
\]

The determination of the rate function \(R_q(\delta)\)
is an asymptotic sphere packing problem in the space \(\zfld_q^n\)
with Hamming distance as \(n\to\infty\), 
which has been considered as a central problem
in combinatorial coding theory (see van Lint \cite{bib:vanLint82a} and Manin \cite{bib:manin81a}).
The asymptotic rate function \(R_q(\delta)\) plays a fundamental role 
in Shannon's information theory as well 
(see \cite{bib:shannon67a} and \cite{bib:shannon67b}).
It is intimately related to the achievable probability of decoding error 
in reliable transmission of information,
which approaches \(0\) exponentially as the number \(n\) of channel uses goes to infinity.
The better a lower bound for the rate function \(R_q(\delta)\) exists,
the larger the error-exponent can be achieved by a sequence of error-correcting codes.

In the recent works \cite{bib:conlon2022a} and \cite{bib:conlon2023a},
Conlon, Fox, {\em et al.} established and explored a tight connection
between set-coloring Ramsey numbers and error-correcting codes.
Therefore, one may expect that the main results obtained in this paper
would be useful for deriving some nontrivial asymptotic bound for set-coloring Ramsey numbers.

It is known that for all integers \(q\geq 2\), 
if \(1-q^{-1}\leq \delta < 1\), then \(R_q(\delta) = 0\). 
The exact value of \(R_q(\delta)\), however, is not known 
for any integer \(q\geq 2\) and real number \(\delta\in (0, 1-q^{-1})\).
A famous lower bound for \(R_q(\delta)\) 
is the following \(q\)-ary asymptotic Gilbert--Varshamov bound 
\cite{bib:gilbert52a}, \cite{bib:varshamov57a} 
\begin{equation}
R_q(\delta) \geq 1 - \delta\log_q(q-1)-\delta\log_q\frac{1}{\delta}-(1-\delta)\log_q\frac{1}{1-\delta}
\dfnt\RGV
\label{eqn:RGVdfnt}
\end{equation}
for all integers \(q\geq 2\) and \(0<\delta<1-q^{-1}\).

Although it was obtained in the 1950s by simple combinatorial arguments in \cite{bib:gilbert52a}
or with the probabilistic method in \cite{bib:varshamov57a},
the \(q\)-ary asymptotic Gilbert--Varshamov bound \(\RGV\) was thought for a long time
that it would be equal to the rate function \(R_q(\delta)\).

Let \(\Fq\) be a finite field with \(q\) elements, where \(q\) is a prime power.
A \(q\)-ary code \(\C\subseteq\Fq^n\) with length \(n\) is called linear 
if it is a subspace of the vector space \(\Fq^n\) over \(\Fq\).
In 1981, Goppa~\cite{bib:goppa81a} found a remarkable connection  
between algebraic curves over \(\Fq\) and \(q\)-ary linear codes
and showed, through the Riemann--Roch theorem, that a better code can be obtained
from an algebraic curve over \(\Fq\) with more rational points given its genus.
In 1982, by using Goppa codes \cite{bib:goppa81a} from algebraic curves over finite fields,
Tsfasman, Vladut, and Zink \cite{bib:tsfasman82a} employed the classic modular curves and the Shimura curves, 
and successfully improved the Gilbert--Varshamov bound \(\RGV\)
for \(q\) being a prime power of \(p^2\) (\(p\geq 7\))
and of \(p^4\) (\(p\geq 3\))
and \(\delta\) in an open subinterval of \((0,1-q^{-1})\).
This discovery was the first breakthrough after over twenty years of searching for a lower bound 
to beat the Gilbert--Varshamov bound \(\RGV\) for certain \(q\) and \(\delta\).

We are going to present briefly theory of Goppa codes from algebraic curves over finite fields.
The following work in this paper on lower bounds of the rate function
will be based on Lenstra codes from geometry of numbers, which can be viewed as
being analogous to Goppa codes.

In 1982, Manin \cite{bib:manin81a} firstly described such connection for algebraic geometers with some arising unsolved problems for drawing their attention.
In 1982, Ihara \cite{bib:ihara81a} defined the following number
\begin{equation}
A(q) = \limsup_{g\rightarrow\infty}\frac{N_q(g)}{g}
\label{eqn:Ihara_A_q}
\end{equation}
where \(N_q(g)\) is the maximum number of rational points 
that a non-singular absolutely irreducible algebraic curve over \(\Fq\) of genus \(g\)
can have.
Ihara's number \(A(q)\) is well-defined by the Hasse--Weil theorem,
namely, the Riemann hypothesis for algebraic curves over finite fields
\cite{bib:stichtenoth09a}.
Ihara \cite{bib:ihara81a} obtained that 
\(A(q)\geq\sqrt{q}-1\) when \(q = p^{2m}\), \(p\) prime and \(m\in\nfld\), 
and provided an upper bound for \(A(q)\) for all \(q\) being a prime power.
It has been proved that Goppa codes imply (see Manin \cite[p.~718]{bib:manin81a})
\begin{equation}
R_q(\delta) \geq 1 - \delta -\frac{1}{A(q)}
\label{eqn:RqAq}
\end{equation}
for all \(q\) being a prime power
and \(0 < \delta < 1\).
According to Manin \cite[p.~719]{bib:manin81a}, 
the work in Tsfasman, Vladut, and Zink \cite{bib:tsfasman82a}, 
independently of the earlier Ihara's work \cite{bib:ihara75a}, \cite{bib:ihara79a}
in the 1970s on towers of curves over \(\Fq\),
showed that \(A(q)\geq\sqrt{q}-1\) when \(q = p^2\) or \(q = p^4\), \(p\) prime, 
and hence from (\ref{eqn:RqAq}) that the Gilbert--Varshamov bound \(\RGV\) was improved
for \(q = p^2\) (\(p\geq 7\))
or \(q = p^4\) (\(p\geq 3\))
and \(\delta\) in an interval of \((\delta_1,\delta_2)\) including \(\frac{q-1}{2q-1}\)
(see \cite[p.~28]{bib:tsfasman82a}).
It is noted that a better upper bound of \(A(q)\) was also obtained,
which is called as the Drinfeld--Vladut bound \cite{bib:vladut83a}, 
asserting that \(A(q)\leq\sqrt{q}-1\) for all~\(q\) being a prime power.
It had improved a previous upper bound of Ihara \cite{bib:ihara81a}. 
The Drinfeld--Vladut bound is tight when \(q = p^{2m}\), \(p\) prime and \(m\in\nfld\), 
from the aforementioned lower bound of \(A(q)\) due to Ihara~\cite{bib:ihara81a}.

Goppa codes from algebraic curves over finite fields are also called as {\em algebraic geometry codes}
(see the monographs \cite{bib:tsfasman07a} and \cite{bib:tsfasman19a}).
It is well-known that algebraic curves over~\(\Fq\) can be equivalently replaced 
in the language of algebraic function fields over~\(\Fq\)~\cite{bib:stichtenoth09a}.
In the case of \(q\) being an even power of a prime,
Garcia and Stichtenoth \cite{bib:garcia95a}
found an alternative approach with the Artin--Schreier tower of algebraic function fields over \(\Fq\)
to attain the Drinfeld--Vladut bound.
When \(q\) is an odd power of a prime {\em but not a prime} with \(q > 125\),
Bassa, Beelen, Garcia, and Stichtenoth \cite{bib:bassa14a}
showed that the Gilbert--Varshamov bound \(\RGV\) can also be improved
by using appropriate infinite towers of algebraic function fields over finite fields~\(\Fq\).

In 1983, Serre \cite{bib:serre83a}, \cite{bib:serre83b} demonstrated that 
there exists a constant \(c>0\) such that 
\(
A(q) \geq \frac{1}{c}\log q
\)
for any integer \(q \geq 2\) being a prime power, 
by using class field theory for algebraic function fields over \(\Fq\)
(see also the recent monograph by Serre \cite{bib:serre2020c}).
This bound is still the best known lower bound for \(A(q)\) in general for \(q\) being a {\em prime}.
By using (\ref{eqn:RqAq}), we have the following Serre's lower bound
\begin{equation}
R_q(\delta) \geq 1 - \delta -\frac{c}{\log q}
\label{eqn:Serre_asymptotic}
\end{equation}
for any integer \(q \geq 2\) being a prime power and \(0<\delta<1-q^{-1}\).
It was shown in Niederreiter and Xing \cite{bib:niederreiter01a} 
that the constant \(c\) in  (\ref{eqn:Serre_asymptotic}) 
can be chosen as \(96\log 2\) for \(q\geq 2\) being a prime power 
or \(75\log 2\) for \(q\geq 2^{24}\) being a prime power.
With this explicit constant \(c\),
it can be checked that Serre's lower bound  (\ref{eqn:Serre_asymptotic})
is  {\em not} better than the Gilbert--Varshamov bound \(\RGV\) for any integer \(q\) being a {\em prime}
and  \(0<\delta<1-q^{-1}\).

In 1985, Zinoviev and Litsyn \cite{bib:zinoviev85a} demonstrated that
the algebraic geometry codes can be modified by an alphabet reduction procedure
to improve the \(q\)-ary Gilbert--Varshamov bound \(\RGV\)
at the specific tangent point
\(\delta=\delta_0 \dfnt (q-1)/(q'+q-1) \in (0,1)\) 
of the curve \(\RGV\) for all \(q\geq 46\),
where \(q'\) is the cardinality of the alphabet for the applied algebraic geometry codes
(see also \cite{bib:litsyn86a}, \cite[p.~267, Theorem~4.5.36.]{bib:tsfasman07a}).
Note that \(q'\) is a prime power {\em but not a prime} in the application \cite{bib:zinoviev85a}.
It is mentioned that a relevant and more recent work is presented by Elkies in 2001 \cite{bib:elkies01a}.

We note that if \(q=p\) is a generic {\em prime},
it is not possible that
the tangent point \(\delta_0\)
can cover the value \(\delta=1/2\) for almost all \(p\).
If \(\delta_0 =1/2\) or is very close to \(1/2\), then the prime \(p = q'+1\) or is very close to \(q'+1\).
Assume that \(p \leq N\), where the positive integer \(N\) is large. 
There are at least \(\Omega(N/\log N) \) primes \(p\) between 2 and \(N\),
while there are at most \(O(N^{1/2}/\log N)\) prime powers \(q'\),
{\em which is not a prime},  between 1 and \(2N\).
Therefore, it cannot be true that  \(p = q'+1\) or is very close to \(q'+1\)
for almost all primes \(p\).
Hence, there are {\em infinitely many} primes \(p\) for which
\(\delta_0\) cannot cover the value \(\delta=1/2\) for \(q=p\).
Next, if \(q\) is a generic {\em non-prime-power} positive integer, then the finite field \(\Fq\)
with \(q\) elements does not exist. 
In the same spirit, we can see that there are {\em infinitely many} such \(q\) for which
the tangent point \(\delta_0\) cannot cover the value \(\delta=1/2\).

Based on the above analysis of previously obtained results in the existing literature,
it can be stated that the \(q\)-ary Gilbert--Varshamov bound \(\RGV\)
at \(\delta=1/2\), i.e., \(R_\mathrm{GV}(1/2,q)\),
remains the currently largest known lower bound of  the rate function \(R_q(\delta)\)
at \(\delta=1/2\), i.e., \(R_q(1/2)\),
for {\em infinitely many} positive integers \(q\) 
which is a generic prime and which is a generic non-prime-power integer.


\subsection{Main results in this paper}

The aim of this paper is to establish the following theorem,
which provides a new lower bound for the rate function \(R_q(\delta)\) with \(0<\delta<1\).

\begin{TT}
There exists an effectively computable absolute constant \(C>0\) such that
\begin{equation}
  R_q(\delta) 
   >  
     1 - \delta - C \frac{\log\log q}{q^{1/6}}
\label{eqn:liang_constant_C}
\end{equation}
for all positive integers \(q\geq 3\) and \(0<\delta<1\).
\label{thrm:thrm_liang_constant_C}
\end{TT}

Now, we can see that the \(q\)-ary Gilbert--Varshamov bound \(\RGV\) with \(\delta\in(0,1)\)
can be improved for all but {\em finitely~many} positive integers \(q\). 
It is clear that the \(q\)-ary asymptotic Gilbert--Varshamov bound given in (\ref{eqn:RGVdfnt})
can be represented by
\begin{equation}
    \RGV = 1 - \delta -\frac{h(\delta)}{\log q} + O\left(\frac{1}{q \log q}\right)
\label{eqn:RGV_asymptotic}
\end{equation}
as \(q\rightarrow\infty\), where 
\(
    h(\delta) = -\delta\log\delta-(1-\delta)\log(1-\delta),
\)
\(0 < \delta < 1\),
is the entropy function. 
One can notice that the above new lower bound given in~(\ref{eqn:liang_constant_C})
for the rate function \(R_q(\delta)\)
has a denominator with \(q^{1/6}\) of {\em algebraic} order in~\(q\),
while the \(q\)-ary Gilbert--Varshamov bound given in~(\ref{eqn:RGV_asymptotic})
has that with \(\log q\)  of only {\em logarithmic} order in~\(q\).
Consequenly,
for an arbitrarily given number \(\ve\in(0,\frac{1}{2})\),
there exists a positive integer \(N(\ve)\) such that
the lower bound given in~(\ref{eqn:liang_constant_C}) for the rate function \(R_q(\delta)\)
exceeds the \(q\)-ary Gilbert--Varshamov bound  \(\RGV\) given in~(\ref{eqn:RGV_asymptotic}),
and hence,  \(R_q(\delta)  > R_\mathrm{GV}(\delta,q)\), for all \(q\geq N(\ve)\) and \(\delta\in[\ve,1-\ve]\).

In particular, we will prove the following {\em explicit} result at \(\delta=1/2\).

\begin{TT}
For all integers \(q > \exp(29)\), we have
\begin{equation}
R_q(1/2)  > R_\mathrm{GV}(1/2,q).
\label{eqn:liang_LB}
\end{equation}
\label{thrm:liang}
\end{TT}

The rate function \(R_q(\delta)\) has the following \(q\)-ary asymptotic Plotkin upper bound 
\cite{bib:plotkin60a},
\cite[p.~68]{bib:vanLint82a}
\[
    R_q(\delta) \leq 1-\delta-\frac{\delta}{q-1}
\]
for all positive integers \(q\geq 2\) and \(0<\delta<1-q^{-1}\),
which is the best known upper bound for large \(q\).
Therefore, we can define the growth on the rate function \(R_q(\delta)\) for \(\delta\in(0,1)\) fixed and \(q\) growing large by
\begin{equation}
\eta(\delta)= \liminf_{q\rightarrow\infty}\frac{1}{\log q}\log[1-\delta-R_q(\delta)]^{-1}.
\label{eqn:eta_delta}
\end{equation}
The \(q\)-ary Plotkin upper bound implies that the growth satisfies \(\eta(\delta)\leq 1\).
On the other hand, 
the \(q\)-ary Gilbert--Varshamov lower bound given in (\ref{eqn:RGV_asymptotic})
only gives the trivial bound \(\eta(\delta)\geq 0\).
From the new lower bound given in (\ref{eqn:liang_constant_C}) for the rate function \(R_q(\delta)\)
in Theorem~\ref{thrm:thrm_liang_constant_C},
we obtain a {\em nontrivial} lower bound for the growth \(\eta(\delta)\) as follows.

\begin{TT}
For every real number \(\delta\in(0,1)\), we have
\begin{equation}
\frac{1}{6} \leq \eta(\delta) \leq 1.
\label{eqn:liang_1over6}
\end{equation}
\label{thrm:liang_1over6}
\end{TT}

Though the {\em exact} value of \(R_q(\delta)\) is still unknown 
for any given integer \(q\geq 2\) and real number \(\delta\in (0, 1-q^{-1})\)
as previously mentioned,
the order of growth on the rate function \(R_q(\delta)\) for \(\delta\in(0,1)\) fixed
and \(q\) growing large
is now determined qualitatively.
The \(q\)-ary Gilbert--Varshamov bound given in (\ref{eqn:RGV_asymptotic}) 
indicates that the reciprocal of the difference between \(1-\delta\) and \(R_q(\delta)\)
has asymptotically at least a {\em logarithmic order} in~\(q\).
The obtained lower bound for the growth \(\eta(\delta)\) given in (\ref{eqn:liang_1over6})
together with the \(q\)-ary Plotkin upper bound
demonstrates that the {\em asymptotic order of growth}
of the reciprocal of the difference between \(1-\delta\) and \(R_q(\delta)\)
in~\(q\) must be {\em algebraic} rather than logarithmic.

The new lower bounds in Theorem~\ref{thrm:thrm_liang_constant_C} and Theorem~\ref{thrm:liang}
are based on theory of Lenstra codes from geometry of numbers,
which can be viewed as an analogue of theory of Goppa codes from algebraic curves over finite fields.
In Section~\ref{sec:Lenstra}, we will present a detailed account of a variant of Lenstra codes.
In Section~\ref{sec:Hilbert}, a theory of infinite towers of Hilbert class fields and
its combination with Lenstra codes are developed.
The theories presented in Section~\ref{sec:Lenstra} 
and in Section~\ref{sec:Hilbert} enable us to obtain an applicable lower bound 
for the rate function \(R_q(\delta)\) in Section~\ref{sec:lowerbound}, 
which will be used in proofs of Theorem~\ref{thrm:thrm_liang_constant_C} and Theorem~\ref{thrm:liang}.
With those preparations,  
we will establish proofs of Theorem~\ref{thrm:thrm_liang_constant_C} and Theorem~\ref{thrm:liang}
in Section~\ref{sec:ProofLiang1over6} and Section~\ref{sec:ProofLiang}, respectively.


\section{Codes from geometry of numbers}
\label{sec:Lenstra}

Let \(K\) be an algebraic number field with degree \(m = [K:\Qfld]\) over the field \(\Qfld\) 
of rational numbers.
If \(K\) has \(s\) real embeddings and \(t\) pairs of complex conjugate embeddings, then we have 
\(s+2t = m\).
We denote by \(O_K\) the ring of integers of \(K\)
and by \(\Delta_K\) the discriminant of \(O_K\)
over the ring \(\zfld\) of rational integers.
In Minkowski's theory of geometry of numbers,
a ring of algebraic integers is viewed as a lattice 
in the \(m\)-dimensional Euclidean space \(\rfld^m\) (see \cite[\S 1.5]{bib:neukirch99a}).

We let \(\sigma_1, \cdots,\sigma_s: K\emb\rfld\) be the \(s\) real embeddings 
and \(\sigma_{s+1}, \ol{\sigma_{s+1}},\cdots,\sigma_{s+t}, \ol{\sigma_{s+t}}: K\emb\cfld\) 
the \(t\) pairs of complex conjugate embeddings of \(K\).
For each integer \(a\in O_K\), let
\(x_i = \sigma_i(a)\) for \(1\leq i\leq s\) when \(s\geq 1\),
and \(y_j + \IU z_j = \sigma_{s+j}(a)\) for \(1\leq j\leq t\) when \(t\geq 1\),
where \(\IU\) denotes the imaginary unit in the field \(\cfld\) of complex numbers.
Then, \(a\in O_K\) can be mapped to a point 
\(x = (x_1,\dotsc,x_m) = (x_1,\dotsc,x_s,y_1,z_1,\dotsc,y_t,z_t)\) in \(\rfld^m\),
denoted by \(x = \Lambda(a): O_K\to\rfld^m\). 
We have \(y_j = x_{s+2j-1}\) and \(z_j = x_{s+2j}\) for \(1\leq j\leq t\).

The image of \(O_K\) under the map \(\Lambda: O_K\to\rfld^m\)
is denoted by \(\Lambda_K\).
Then, \(\Lambda_K\) is a lattice in \(\rfld^m\).
It means that \(\Lambda_K\) is an additive subgroup of \((\rfld^m,+)\)
and it is discrete in the sense that
for any bounded subset \(U\subset\rfld^m\),
the intersection \(U\cap\Lambda_K\) is a finite subset of \(\rfld^m\).
Let \(T\) be a fundamental domain of the lattice \(\Lambda_K\).
Then, its volume with Lebesgue measure is given by
\(
   \vol(T) = 2^{-t}\sqrt{|\Delta_K|}
\)
(see, e.g., \cite[Theorem~9.4]{bib:stewart02a}).

For a finite set \(\Omega\), its cardinality is denoted by \(|\Omega|\).
The following lemma is an averaging argument in geometry of numbers
\cite{bib:lenstra86a}, \cite{bib:siegel45a}.

\begin{LL} {\em (Lenstra--Siegel)}
If \(U\) is a nonempty bounded subset of \(\rfld^m\) which is Lebesgue measurable with
volume denoted by \(\vol(U)\),
then there exists \(\tau\in T\) such that
\[
|(\tau+U)\cap\Lambda_K|\geq\frac{\vol(U)}{\vol(T)}.
\]
\label{lemma:Lenstra_Siegel}
\end{LL}

{\em Proof:}
Since \(\wt{T}\dfnt -T\) is also a fundamental domain of the lattice \(\Lambda_K\),
by \cite[Lemma~6.2]{bib:stewart02a} we know that
each point \(z\in\rfld^m\) lies in exactly one of the cosets \(\tau+\Lambda_K\)
where \(\tau\in \wt{T}\).

Let \(\chi_U(z)\) for \(z\in\rfld^m\) be the characteristic function of the Lebesgue measurable
subset \(U\subseteq\rfld^m\)
defined as \(1\) if \(z\in U\) and \(0\) otherwise.

We evaluate the following integral
\begin{eqnarray}
  \int_{\tau\in \wt{T}} |(-\tau+U)\cap\Lambda_K| d\tau 
  & = & 
  \int_{\tau\in \wt{T}}\sum_{y\in\Lambda_K} |(-\tau+U)\cap\{y\}|d\tau
   \nonumber \\
   & = & 
  \int_{\tau\in \wt{T}}\sum_{y\in\Lambda_K}  \chi_U(\tau+y) d\tau
   \nonumber \\
   & = & 
  \int_{z\in \rfld^m}  \chi_U(z) dz
   \nonumber \\
  & = & \vol(U). \nonumber 
\label{eqn:integral_argument}
\end{eqnarray}
One can choose any \(\tau\in \wt{T} = -T\), i.e., \(-\tau\in T\)
such that the cardinality
\(|(-\tau+U)\cap\Lambda_K|\) is maximal, which is at least \(\vol(U)/\vol(T)\).
The proof of this lemma is complete.
\eof

In the following, we introduce a variant of Lenstra's construction of 
codes from geometry of numbers \cite{bib:lenstra86a}.
We assume that \(q\) and \(r\) are two arbitrarily given positive integers
satisfying \( 2\leq r \leq q\).

We choose \(U\) in Lemma~\ref{lemma:Lenstra_Siegel} as the set of \(x\in\rfld^m\)
with
\[
  0  < x_i < \rho\dfnt 2^{-t/m}r^{G/m}
\]
for \(1\leq i\leq m\), denoted by \(U_G\),
where \(G\in\nfld\) is a parameter. 
It is clear that \(U_G\) is an open subset of \(\rfld^m\) with
\(
\vol(U_G) = \rho^m = 2^{-t}r^G.
\)
By Lemma~\ref{lemma:Lenstra_Siegel}, there exists \(\tau\in T\) such that
\begin{equation}
|(\tau+U_G)\cap\Lambda_K| \geq \frac{r^G}{\sqrt{|\Delta_K|}}.
\label{eqn:tau_U}
\end{equation}
We define the following finite subset of \(O_K\)
\[
      \Omega_G = \{ a\in O_K \mid \Lambda(a) \in (\tau+U_G)\cap\Lambda_K \}
\]
whose cardinality is equal to \(|(\tau+U_G)\cap\Lambda_K|\).
By (\ref{eqn:tau_U}), we have \(|\Omega_G|\geq r^G/\sqrt{|\Delta_K|}\).

Let \(P\) be a nonzero prime ideal of \(O_K\).
Then, the quotient ring \(O_K/P\) is a finite field, which is called the residue class field
of \(O_K\) modulo \(P\).
The norm of \(P\) is defined by \(\norm(P) = |O_K/P|\).
We define the number
\begin{equation}
   N_{r,q}(K) = |\{P\subseteq O_K \mid r \leq \norm(P) \leq q \}|
\label{eqn:N_r_q_K}
\end{equation}
where \(P\) runs through the nonzero prime ideals of \(O_K\). 
Note that \(N_{r,q}(K)\) is the cardinality of the given subset which has finitely many elements.

In the case that \( n = N_{r,q}(K) \geq 1\), 
we let \(P_1,\cdots,P_n\) denote the distinct nonzero prime ideals of \(O_K\)
whose norms are between \(r\) and \(q\).
It is ready to define the following map \(\psi:\Omega_G\to\zfld_q^n\)
\[
    \psi(a) = (a(P_1),\dotsc,a(P_n)), \qquad \mbox{for all } a\in\Omega_G
\]
where \(a(P_i) = a + P_i\in O_K/P_i\)  is the residue class of \(a\in O_K\) 
modulo \(P_i\) for \(1\leq i\leq n\).
By choosing an arbitrary injective map \(O_K/P_i\to\zfld_q\) for each \(i\) with \(1\leq i\leq n\), 
we understand that \(\psi(a)\in\zfld_q^n\).

The image of the map \(\psi:\Omega_G\to\zfld_q^n\),
denoted by \(\code(P_1,\cdots,P_n;G)\),
is a subset of \(\zfld_q^n\),
which is called a \(q\)-ary code from geometry of numbers
or a \(q\)-ary Lenstra code. 
It has some properties as described below.

\begin{PP} 
Let \(K\) be an algebraic number field and \(r,q\in\nfld\) satisfying \(2\leq r\leq q\).
Assume \( n = N_{r,q}(K) \geq 1\). 
If \(r^n \geq \sqrt{|\Delta_K|}\),
then for any \(G\in\nfld\) with 
\(r^G \geq \sqrt{|\Delta_K|}\) and \(G\leq n\),
the Lenstra code \(\code(P_1,\cdots,P_n;G)\)
is a \(q\)-ary \((n,M,d)\) code with 
\(M\geq r^G/\sqrt{|\Delta_K|}\)
and
\(d\geq n+1-G\).
\label{prop:Lenstra_code}
\end{PP}

{\em Proof:}
Let \(a,b\in\Omega_G\subseteq O_K\) with \(a\neq b\), whose corresponding codewords are
\(c_a = \psi(a)\in\zfld_q^n\) and \(c_b = \psi(b)\in\zfld_q^n\), respectively.

Let \(w = \dH(c_a,c_b)\) be the Hamming distance between \(c_a\) and \(c_b\).
Then, \(w = n-\ell\), where \(\ell\) is the number of \(P_i\)'s, \(1\leq i\leq n\), 
such that \(a(P_i) = a + P_i = b(P_i) = b + P_i\).
Without loss of generality, we assume \(P_1,\cdots,P_\ell\) being those nonzero prime ideals.
Hence, \(a-b\in P_i\) for \(1\leq i\leq\ell\).

Now, each \(P_i\), \(1\leq i\leq\ell\), is a prime factor of the prime ideal factorization 
in~\(O_K\) of the nonzero principal ideal \((a-b)O_K \subseteq O_K\).
As a consequence, the norm of \((a-b)O_K\),
defined as the cardinality of the quotient ring \(O_K/((a-b)O_K)\),
is at least 
\(
   \prod_{i=1}^\ell\norm(P_i).
\)

Note that the norm of the nonzero principal ideal
\((a-b)O_K\) is equal to the absolute value of the field norm
\(\NKQ(a-b)\) \cite[Corollary~5.10]{bib:stewart02a}.
Thus, we have
\begin{equation}
| \NKQ(a-b) | 
\geq \prod_{i=1}^\ell\norm(P_i)\geq r^\ell.
\label{eqn:norm_lower_bound}
\end{equation}

Let \(x = \Lambda(a)\) with \(x = (x_1,\dotsc,x_m) = (x_1,\dotsc,x_s,y_1,z_1,\dotsc,y_t,z_t)\)
and \(\hat{x} = \Lambda(b)\) with \( \hat{x} = (\hat{x}_1,\dotsc,\hat{x}_m) = (\hat{x}_1,\dotsc,\hat{x}_s,\hat{y}_1,\hat{z}_1,\dotsc,\hat{y}_t,\hat{z}_t)\)
under the map \(\Lambda: O_K\to\rfld^m\).

Since \(a,b\in\Omega_G\), we have \(x,\hat{x}\in(\tau+U_G)\cap\Lambda_K\), 
where \(\tau = (\tau_1,\dotsc,\tau_m)\in T\) is given in (\ref{eqn:tau_U}).
Therefore, \(x-\tau,\hat{x}-\tau\in U_G\),
namely,
\(
 0 < x_i - \tau_i < \rho
\)
and
\(
 0 < \hat{x}_i - \tau_i < \rho
\)
for \(1\leq i\leq m\).
Hence, 
\(
 |x_i - \hat{x}_i| = |(x_i - \tau_i)- (\hat{x}_i - \tau_i)| < \rho
\)
for \(1\leq i\leq m\)
and
\(
   |y_j - \hat{y}_j|^2 + |z_j - \hat{z}_j|^2
 = |x_{s+2j-1} - \hat{x}_{s+2j-1}|^2 + |x_{s+2j} - \hat{x}_{s+2j}|^2
 < 2 \rho^2
\)
for \(1\leq j\leq t\).

We estimate the absolute value of the norm \(\NKQ(a-b)\) as follows
\begin{eqnarray}
  | \NKQ(a-b) | 
  & = & 
  \prod_{i=1}^s|\sigma_i(a-b)| \times \prod_{j=1}^t|\sigma_{s+j}(a-b)|^2
   \nonumber \\
   & = & 
  \prod_{i=1}^s|\sigma_i(a)-\sigma_i(b)| \times \prod_{j=1}^t|\sigma_{s+j}(a)-\sigma_{s+j}(b)|^2
   \nonumber \\
   & = & 
  \prod_{i=1}^s|x_i - \hat{x}_i| \times 
                    \prod_{j=1}^t |(y_j + \IU z_j) - (\hat{y}_j + \IU \hat{z}_j)|^2
   \nonumber \\
   & = & 
  \prod_{i=1}^s|x_i - \hat{x}_i| \times  \prod_{j=1}^t[|y_j - \hat{y}_j|^2 + |z_j - \hat{z}_j|^2]
   \nonumber \\
  & < & \rho^s \times (2\rho^2)^t = r^G. 
\label{eqn:absolute_value_norm}
\end{eqnarray}
Then, it follows from (\ref{eqn:norm_lower_bound}) and (\ref{eqn:absolute_value_norm})
that \(r^G > r^\ell\), i.e.,
\(G\geq \ell+1\).
Hence, \(w = n-\ell\geq n+1-G\). 
Since \(G\leq n\), we have \(\dH(c_a,c_b) = w \geq 1\),
which implies that \(c_a\neq c_b\).
Thus, the map \(\psi:\Omega_G\to\zfld_q^n\) is injective.

Consequently, 
the code \(\code(P_1,\cdots,P_n;G)\),
defined as the image of the map \(\psi:\Omega_G\to\zfld_q^n\),
has the size of
\(M = |\code(P_1,\cdots,P_n;G)| = |\Omega_G|\geq r^G/\sqrt{|\Delta_K|}\)
and the minimum Hamming distance of
\(d\geq n+1-G\).
Therefore, the proof is complete.
\eof

The code \(\code(P_1,\cdots,P_n;G)\) in Proposition~\ref{prop:Lenstra_code} 
is actually a {\em punctured} version \cite[p.~28]{bib:macwilliams98a} 
of the code from geometry of numbers originally introduced by Lenstra~\cite{bib:lenstra86a} in 1986.
The length of the code constructed in \cite{bib:lenstra86a} is equal to
\begin{equation}
   \hat{N}_{r,q}(K) = s+t+|\{P\subseteq O_K \mid r \leq \norm(P)^{k(P)} \leq q \; \mbox{for some } 
     k(P)\in\nfld\}|
\label{eqn:hat_N_r_q}
\end{equation}
where \(P\) runs through the nonzero prime ideals of \(O_K\). 
It is clear that
\(\hat{N}_{r,q}(K) \geq s+t+N_{r,q}(K) \geq N_{r,q}(K)\).
The Lenstra code in \cite{bib:lenstra86a} has involved all real and complex embeddings of 
the number field \(K\) and possibly more nonzero prime ideals \(P\subseteq O_K\)
in the construction. 
The code \(\code(P_1,\cdots,P_n;G)\) in Proposition~\ref{prop:Lenstra_code} 
is constructed via {\em puncturing} the original Lenstra code by
keeping only coordinates of each codeword 
at prime ideals with their norms between \(r\) and \(q\).
The {\em punctured} Lenstra code is employed because in this case 
theory of codes from geometry of numbers is simpler,
and it is adequate for our purpose of 
establishing Theorem~\ref{thrm:liang} and Theorem~\ref{thrm:liang_1over6}.

It is noticed that Lenstra \cite{bib:lenstra86a} only gave a complete proof of the assertion
regarding properties of \((n,M,d)\) for codes from geometry of numbers 
{\em in the case} of \(K\) with \(s=[K:\Qfld]\) and \(t=0\), namely,  \(K\) being a totally real field.
Nonetheless, we have given a proof of Proposition~\ref{prop:Lenstra_code}
as to properties of the simpler {\em punctured} Lenstra codes
for an {\em arbitrary} number field \(K\) with \(t=0\) or \(t>0\).
Both cases of \(t=0\) and \(t>0\) will be needed for the code constructions in the paper.

It is known that the discriminant \(\Delta_K\) of a number field \(K\) 
is always a nonzero rational integer.
We introduce the number
\[
   N_{r,q}(\Delta) = \max_{K:|\Delta_K| = \Delta}N_{r,q}(K)
\]
for any given \(\Delta\in\nfld\), which is well-defined by the Hermite theorem
\cite[Theorem~3.2.16]{bib:neukirch99a}
asserting that there are only finitely many algebraic number fields \(K\),
up to isomorphism,
with \(|\Delta_K| \leq \Delta\).
Then, for \(r,q\in\nfld\) with \(2\leq r\leq q\), we can define the quantity
\begin{equation}
     A(r,q) = \limsup_{\Delta\rightarrow\infty}\frac{N_{r,q}(\Delta)}{\log\sqrt{\Delta}}.
\label{eqn:Lenstra_A_r_q}
\end{equation}

It is clear that for \(r_1,r_2,q\in\nfld\) with \(2\leq r_1\leq r_2\leq q\), we have
\(A(r_1,q) \geq A(r_2,q)\).
In general, \(A(r,q)\) has the following upper bound.

\begin{LL} 
For all \(r,q\in\nfld\) with \(2\leq r\leq q\), we have
\[
   0 \leq A(r,q) < 7\frac{q}{\log q}.
\]
\label{lemma:A_r_q}
\end{LL}

{\em Proof:}
It is obvious that \(A(r,q) \geq 0\). We can assume that \(A(r,q) > 0\).

Let \(K\) be an algebraic number field.
For each prime ideal \(P\) in the subset
given in~(\ref{eqn:N_r_q_K}) defining the number \(N_{r,q}(K)\), 
it lies over some prime number \(p\) in the sense that
\(P\cap \zfld = p\zfld\).
Then, \(\norm(P) = p^f\), where \(f\geq 1\) is the inertia degree of \(P\) over \(p\zfld\).
Since \(\norm(P)\leq q\), we have \(p\leq q\). 
On the other hand, by the fundamental identity in extension of algebraic number fields
\cite[Proposition 1.8.2]{bib:neukirch99a},
there are at most \([K:\Qfld]\) prime ideals of \(O_K\) lying over \(p\zfld\). 
Therefore, we obtain
\begin{equation}
     N_{r,q}(K) \leq [K:\Qfld]\pi(q)
\label{eqn:N_r_q_pi_q}
\end{equation}
where \(\pi(q)\) is the prime-counting function.

We need Minkowski's bound on the discriminant \cite[Proposition 3.2.14]{bib:neukirch99a}
\begin{equation}
    \sqrt{|\Delta_K|} \geq \frac{m^m}{m!}\left(\frac{\pi}{4}\right)^{m/2}
\label{eqn:Minkowski_bound}
\end{equation}
where \(m = [K:\Qfld]\).

Let \(K_i\), \(i\in\nfld\), be a sequence of number fields such that
\(\Delta_i = |\Delta_{K_i}| > 1\)
and
\(N_{r,q}(\Delta_i) = N_{r,q}(K_i)\),
which satisfy
\[
    \lim_{i\to\infty}  \Delta_i = \infty       \mbox{ and }
    \lim_{i\to\infty}\frac{N_{r,q}(K_i)}{\log\sqrt{\Delta_i}}
    = A(r,q) > 0.
\]
This implies \(N_{r,q}(K_i)\to\infty\), and hence, from (\ref{eqn:N_r_q_pi_q}), 
that \(m_i = [K_i:\Qfld]\to\infty\) as \(i\to\infty\).
Therefore, by using Minkowski's bound (\ref{eqn:Minkowski_bound})
and Stirling's formula for factorials, we get
\begin{equation}
    A(r,q) \leq \frac{1}{1-\log\frac{2}{\sqrt{\pi}}}\pi(q)=1.1373\cdots\pi(q).
\label{eqn:A_r_q_upper_bound}
\end{equation}
Since \(\pi(q) < 6\frac{q}{\log q}\) for all \( q\geq 2\) \cite[Theorem~4.6]{bib:apostol76a}, 
the proof of this lemma is finished.
\eof

\begin{LL} 
For all positive integers \( q\geq 2\), we have
\[
   0 \leq A(q,q) \leq \frac{1}{1-\log\frac{2}{\sqrt{\pi}}}=1.1373\cdots.
\]
\label{lemma:A_q_q_bigO}
\end{LL}

{\em Proof:}
We employ the above proof of Lemma~\ref{lemma:A_r_q}. 
In the specific case of \(r=q\), 
there is either no any or only one prime number \(p\) such that
\(\norm(P) = p^f = q\) for the given \(q\).
Therefore, we obtain
\[
     N_{r,q}(K) \leq [K:\Qfld]
\]
instead of (\ref{eqn:N_r_q_pi_q}). 
Then, the upper bound in this lemma follows, instead of (\ref{eqn:A_r_q_upper_bound}).
\eof

With the quantity \(A(r,q)\) defined in (\ref{eqn:Lenstra_A_r_q}),
the Lenstra code with properties as detailed in Proposition~\ref{prop:Lenstra_code} 
can provide the following lower bound for the rate function \(R_q(\delta)\).

\begin{PP} 
For all \(r,q\in\nfld\) with \(2\leq r\leq q\), we have
\begin{equation}
   R_q(\delta) \geq (1-\delta)\frac{\log r}{\log q} - \frac{1}{A(r,q)\log q}
\label{eqn:Rq_A_r_q_lower_bound}
\end{equation}
for \(0 < \delta < 1\).
\label{prop:Rq_A_r_q}
\end{PP}

{\em Proof:}
We can assume that 
\[
     0 < \delta < 1 - \frac{1}{A(r,q)\log r}
\]
for otherwise the right-hand side of (\ref{eqn:Rq_A_r_q_lower_bound})
is less than or equal to zero.

Let \(K_i\), \(i\in\nfld\), be a sequence of number fields such that
\(\Delta_i = |\Delta_{K_i}| > 1\)
and
\(N_{r,q}(\Delta_i) = N_{r,q}(K_i)\)
satisfying
\begin{equation}
    \lim_{i\to\infty}  \Delta_i = \infty       \mbox{ and }
    \lim_{i\to\infty}\frac{N_{r,q}(K_i)}{\log\sqrt{\Delta_i}}
    = A(r,q) > \frac{1}{\log r}.
\label{eqn:n_i_Lenstra}
\end{equation}

In what follows, we consider for \(i\in\nfld\) sufficiently large.

Let \(n_i = N_{r,q}(K_i)\). Then, from (\ref{eqn:n_i_Lenstra}), 
\(n_i\to\infty\) as \(i\to\infty\).
Moreover, \(r^{n_i} > \sqrt{|\Delta_{K_i}|}\).

We choose \(G_i = \lceil (1-\delta)n_i \rceil \in\nfld\). Then, \(G_i < n_i\).
It is clear that
\[
       \lim_{i\to\infty} \frac{G_i}{\log\sqrt{|\Delta_{K_i}|}}
    =  \lim_{i\to\infty} \frac{G_i}{n_i}\frac{N_{r,q}(K_i)}{\log\sqrt{\Delta_i}}
    =  (1-\delta)A(r,q)
    >  \frac{1}{\log r}.
\]
Therefore, \(r^{G_i} > \sqrt{|\Delta_{K_i}|}\) for \(i\in\nfld\) sufficiently large.

By Proposition~\ref{prop:Lenstra_code},
the Lenstra code \(\code(P_1,\cdots,P_{n_i};G_i)\) is a \(q\)-ary \((n_i,M_i,d_i)\) code with 
\[
   M_i \geq r^{G_i}/\sqrt{|\Delta_{K_i}|}
\]
and
\[
   d_i \geq n_i + 1 - G_i \geq \delta n_i
\]
for \(i\in\nfld\) sufficiently large.
Therefore, we have \(\alpha_q(n_i,\delta n_i)\geq r^{G_i}/\sqrt{|\Delta_{K_i}|}\),
and hence
\begin{eqnarray*}
        R_q(\delta) 
   & = &    \limsup_{n\rightarrow\infty}\frac{1}{n}\log_q\alpha_q(n,\delta n)
   \nonumber \\
   & \geq & \limsup_{i\to\infty} \frac{1}{n_i}\log_q\alpha_q(n_i,\delta n_i)
   \nonumber \\
   & \geq &   \lim_{i\to\infty} \frac{1}{\log q} 
        \left[ \frac{G_i}{n_i}\log r - \frac{\log\sqrt{|\Delta_{K_i}|}}{n_i}  \right]
        \nonumber \\
      \nonumber \\
   & = & (1-\delta)\frac{\log r}{\log q} - \frac{1}{A(r,q)\log q}.
\end{eqnarray*}
The proof of this proposition is complete.
\eof

By choosing \(\delta = 1 - q^{-1}\) for which \(R_q(\delta)=0\),
Proposition~\ref{prop:Rq_A_r_q} can yield the following
upper bound of \(A(r,q)\), which depends on both \(r\) and \(q\).

\begin{CC}
For all \(r,q\in\nfld\) with \(2\leq r\leq q\), we have
\[
   0 \leq A(r,q) \leq \frac{q}{\log r}.
\]
\label{coro:Rq_A_r_q_upper_bound}
\end{CC}

It is interesting that this new upper bound of \(A(r,q)\)
cannot be subsumed by (\ref{eqn:A_r_q_upper_bound}),
since it is better than (\ref{eqn:A_r_q_upper_bound}) in the case of \(r=\lfloor q^{0.99} \rfloor\)
and \(q\) being sufficiently large, as seen by the prime number theorem
\cite[p.~74]{bib:apostol76a}.

From (\ref{eqn:Rq_A_r_q_lower_bound}) in Proposition~\ref{prop:Rq_A_r_q},
we observe that the Lenstra code would provide a lower bound on 
the rate function \(R_q(\delta)\), 
if a lower bound on \(A(r,q)\), which is defined by~(\ref{eqn:Lenstra_A_r_q}),
were found.

We are going to make use of infinite towers of Hilbert class fields
in class field theory for algebraic number fields
to achieve a lower bound for the quantity \(A(r,q)\).


\section{Hilbert class field tower and its concatenation with Lenstra codes}
\label{sec:Hilbert}

It was an open question in algebraic number theory, 
going back to Furtw{\"a}ngler and Hilbert,
whether every number field can be embedded in a number field with
class number being \(1\) (i.e., equal to its Hilbert class field).
In 1964, Golod and Shafarevich~\cite{bib:golod64a} provided a negative answer
to this question by showing that the Hilbert class field tower of a number field \(K\)
may be {\em infinite}.
An explicit example of such number field given in \cite[p.~270]{bib:golod64a} is
a quadratic field
\(
   K = \Qfld(\sqrt{-3\cdot 5\cdot 7\cdot 11\cdot 13\cdot 17\cdot 19})
     = \Qfld(\sqrt{-4849845})
\).

Serre \cite{bib:serre83a}, \cite{bib:serre83b} firstly presented a surprising application
of an infinite Hilbert class field tower, with prescribed places splitting completely,
for algebraic {\em function fields} of curves over finite fields
to obtain a general lower bound on Ihara's number
\(A(q)\) defined by~(\ref{eqn:Ihara_A_q})
(see also \cite{bib:niederreiter01a} and \cite{bib:schoof90a} for expositions).

In this section, we follow the idea of Serre to show that one can apply
an infinite tower of Hilbert class fields, with prescribed places splitting completely,
for algebraic {\em number fields}
to obtain a lower bound for the quantity \(A(r,q)\) defined by~(\ref{eqn:Lenstra_A_r_q}).
Then, a lower bound for \(A(r,q)\) attainable from infinite Hilbert class field towers
can yield, through Lenstra codes, a lower bound for the rate function \(R_q(\delta)\).


\subsection{Hilbert class field tower}

We give some preliminaries about infinite towers of Hilbert class fields
of an algebraic number field in class field theory.

A {\em place} (or {\em prime}) \(P\) of an algebraic number field \(K\)
is a class of equivalent valuations, or a normalized valuation, of \(K\).
The archimedean equivalence classes are called {\em infinite} places,
and the non-archimedean ones called {\em finite} places.
Let \(\Pfld_K\) denote the set of places of \(K\),
\(S_{\infty}\) the set of infinite, or {\em archimedean}, places of \(K\),
and \(S_f\) the set of finite, or {\em non-archimedean}, places of \(K\).
Thus, \(\Pfld_K = S_{\infty}\cup S_f\). 
An infinite place \(P\in\Pfld_K\) is referred to by the notation \(P\mid\infty\)
and a finite one by \(P\nmid\infty\).
It is known that the infinite places \(P\mid\infty\) of \(K\)
are obtained from the real and complex embeddings of \(K\).
There are \(s\) real places associated with the \(s\) real embeddings of \(K\)
and \(t\) complex places with the \(t\) pairs of complex conjugate embeddings of \(K\).
Thus, the number of infinite places of \(K\) is given by
\(|S_{\infty}| = s + t\).
The finite places \(P\nmid\infty\) of \(K\) can stand for the nonzero prime ideals \(P\) of \(O_K\)
by the normalized \(P\)-adic exponential valuation
(see, e.g., \cite{bib:neukirch99a}).

Let \(I_K\) denote the multiplicative group of nonzero fractional ideals of \(K\),
and \(P_K\) its subgroup of nonzero principal fractional ideals of \(K\).
The quotient group \(I_K/P_K\), denoted by \(Cl_K\),
is called the ideal class group of \(K\),
which is a {\em finite} abelian group.
The order of the ideal class group \(Cl_K\) is called the class number of \(K\),
denoted by \(h_K\).
Let \(U_K\) be the multiplicative group of units of the ring of integers \(O_K\),
which is a subgroup of the multiplicative group \(K^* = K \setminus \{0\}\). 
It is noteworthy to remark that for the two multiplicative groups of \(K^*\) and \(I_K\),
there is a group homomorphism \(f: K^*\to I_K\) sending \(a\in K^*\) to 
the principal fractional ideal \(aO_K\in I_K\).
Then, \(U_K\) is the kernel of \(f\) denoted by \(\Ker f\),
\(P_K\) the image of \(f\) by \(\Image f\),
and \(Cl_K\) the cokernel of \(f\) by \(\Coker f\).

A main result of class field theory is that,
for an algebraic number field \(K\) of finite degree,
there exists a maximal unramified abelian extension of \(K\), denoted by \(H_K\),
in which the real places of \(K\) stay real.
The number field \(H_K\) is called the {\em Hilbert class field} of \(K\)
(see \cite[p.~399]{bib:neukirch99a}).
The most important property of \(H_K\) is that
the Galois group of the extension \(H_K/K\) is isomorphic to the ideal class group 
\(Cl_K\), i.e., \(\Gal(H_K/K)\iso Cl_K\). There is no finite place  \(P\nmid\infty\)
of \(K\) ramified in \(H_K\). Moreover, if \(\sigma: K\emb\rfld\) is a real embedding of \(K\),
then all extensions of \(\sigma\) to \(H_K\) remain real.
It is clear that \(H_K = K\) if and only if \(h_K = 1\) (i.e., \(O_K\) is a principal ideal domain).

We are ready to define a sequence of field extensions as follows
\[
      K=K_1 \subseteq K_2 \subseteq K_3 \subseteq\cdots
\]
where \(K_{i+1} = H_{K_i}\), the Hilbert class field of \(K_i\), for \(i\in\nfld\).
This sequence of fields is called the {\em Hilbert class field tower} of \(K\).
We say that \(K\) has an {\em infinite} Hilbert class field tower
if \([K_i:K]\to\infty\) as \(i\to\infty\) and {\em finite} otherwise.

We need a notion of \(p\)-rank in group theory \cite{bib:niederreiter01a}.
If \(G\) is a finitely generated abelian group (written {\em additively}) and \(p\) is a prime, 
the \(p\)-rank of \(G\), denoted by \(d_p G\),  
is defined as the dimension of the \(\FF_p\)-vector space \(G/pG\).
If \(G\) is a finite abelian group, then the \(p\)-rank \(d_p G\) is equal to the number of summands 
in the decomposition of the \(p\)-Sylow subgroup of \(G\) as a direct sum of cyclic subgroups.
In general, a finitely generated abelian group \(G\) can be decomposed as a direct sum of
\(G = G_t \oplus G_f\),
where \(G_t\) is the torsion subgroup of finite order of \(G\),
and \(G_f\) is a free abelian subgroup of \(G\).
Then, the \(p\)-rank \(d_p G\) is the sum of the \(p\)-rank \(d_p G_t\)  
plus the rank of the free abelian group \(G_f\).

The Golod--Shafarevich condition for the infinite Hilbert class field tower
of \(K\) (see Roquette \cite[p.~235]{bib:roquette67a}) asserts that if, for some prime \(p\),
\begin{equation}
d_p Cl_K \geq 2 + 2 \sqrt{d_p U_K + 1}
\label{eqn:GS_condition_Cl_U}
\end{equation}
then \(K\) has an infinite Hilbert class field tower.

By Dirichlet's unit theorem \cite[p.~300]{bib:stewart02a}, 
\(U_K\) is a finitely generated abelian group and it follows that
\begin{equation}
    d_p U_K = s+t-1 +\delta_p(K) = |S_\infty|-1 +\delta_p(K)
\label{eqn:Dirichlet_unit}
\end{equation}
where \(\delta_p(K)\) equals \(1\) if \(K\) contains a primitive \(p\)th root of unity 
and \(0\) otherwise.

It is indicated by (\ref{eqn:GS_condition_Cl_U}) that 
the ideal class group \(Cl_K\) and the group of units \(U_K\) of \(K\)
can jointly determine a situation to guarantee
infiniteness of the Hilbert class field tower of \(K\). 
By virtue of the above mentioned homomorphism \(f: K^*\to I_K\),
the Golod--Shafarevich condition (\ref{eqn:GS_condition_Cl_U}) can be rewritten as
\[
    d_p \Coker f \geq 2 + 2 \sqrt{d_p \Ker f + 1}.
\]

The above theory of infinite Hilbert class field towers
can be generalized to including a requirement of
some prescribed finite places {\em splitting completely} in the class field.
The generalized theory can then be connected with Lenstra codes to provide
a new lower bound for the rate function \(R_q(\delta)\).
We need to use the concepts of ad\`{e}le and id\`{e}le in the generalized theory.
The two groups of \(Cl_K\) and \(U_K\) 
in the Golod--Shafarevich condition~(\ref{eqn:GS_condition_Cl_U}) 
will be replaced by a generalized version of theirs.


\subsection{\(S_c\)-Hilbert class field tower}

Let \(K\) be an algebraic number field
and \(S_c\) a finite subset of {\em finite} places of \(K\),
i.e., \(S_c\subseteq S_f\subseteq \Pfld_K\) with \(|S_c|<\infty\).
In the following, we present a route for theory of \(S_c\)-Hilbert class field of \(K\)
(see Tate \cite{bib:tate67a}, Cassels \cite{bib:cassels67a}, 
and Neukirch \cite[Chapter~VI]{bib:neukirch99a}).

For each \(P\in\Pfld_K\), let \(K_P\) denote the completion of the global field \(K\)
and \(K_P^* = K_P \setminus \{0\}\).
If \(P\mid\infty\) is a real or complex place, then
\(K_P = \rfld\) or \(\cfld\), respectively. 
If \(P\nmid\infty\), 
let~\(v_P\) be the normalized discrete valuation of \(K_P\)
with \(v_P(K_P^*) = \Zfld\).
For each \(P\nmid\infty\),
the valuation ring of \(K_P\) is denoted by 
\(O_P = \{a\in K_P\mid v_P(a)\geq 0\}\)
and the group of units of \(O_P\) by
\(U_P = \{a\in K_P\mid v_P(a) = 0\}\).

An ad\`{e}le of \(K\) is an element 
\(\alpha=(\alpha_P)\) of the direct product \(\prod_{P\in\Pfld_K}K_P\)
such that the \(P\)-component \(\alpha_P\in O_P\) 
for all but finitely many places including infinite ones.
The set of ad\`{e}les of \(K\), denoted by \(\Ad_K\), 
is a ring under addition and multiplication defined componentwise.
The number field \(K\) can be embedded into \(\Ad_K\) 
by the diagonal embedding \(K \hookrightarrow \Ad_K\)
sending \(a\in K\) to \(\alpha=(\alpha_P)\in \Ad_K\) with 
\(\alpha_P = a\) for all  \(P\in\Pfld_K\).

The id\`{e}le group of \(K\), denoted by \(J_K\), is the multiplicative group of units of \(\Ad_K\).
An id\`{e}le of \(K\) is an element 
\(\alpha=(\alpha_P)\) of the direct product \(\prod_{P\in\Pfld_K}K_P^*\)
such that the \(P\)-component \(\alpha_P\in U_P\) 
for all but finitely many places including infinite ones.
The multiplicative group \(K^*\) can be embedded into the abelian group \(J_K\) 
by the diagonal embedding \(K^* \hookrightarrow J_K\)
sending \(a\in K^*\) to \(\alpha=(\alpha_P)\in J_K\) with 
\(\alpha_P = a\) for all  \(P\in\Pfld_K\).
Each element of \(K^*\), which is viewed as a subgroup of \(J_K\),
is called a {\em principal id\`{e}le} of~\(K\).
The quotient group \(C_K = J_K/K^*\) is called the {\em id\`{e}le class group} of \(K\).

Before we move forward, 
it is mentioned that the following three group isomorphism theorems \cite[p.~44]{bib:hungerford74a}
are useful in our development of theory.

Let \(G\) be an abelian group (written {\em multiplicatively}).
Then, each subgroup of \(G\) is a {\em normal} subgroup.
Let \(N\) and \(K\) be subgroups of the abelian group \(G\).
We denote by~\(NK\) the subset 
\(\{ab \mid a\in N \mbox{ and } b\in K\}\subseteq G\).
Then, \(NK\) is a subgroup of \(G\) and it includes \(N\) and~\(K\) as its subgroups.
Let \(H\) be a group. 
The {\em first} group isomorphism theorem
asserts that if \(\phi:G\to H\) is a group homomorphism, then
\(G/\Ker\phi\iso\Image\phi\).
The {\em second} group isomorphism theorem claims that 
\(NK/N \iso K/(N\cap K)\).
The {\em third} group isomorphism theorem says that 
if \(N\) is a subgroup of \(K\), then \(K/N\) is a normal subgroup of \(G/N\) and 
\((G/N)/(K/N) \iso G/K\).

Note that each finite place \(P\in S_f\) can be identified by
a nonzero prime ideal \(P\) of~\(O_K\).
There is a surjective homomorphism \(\pi_f: J_K\to I_K\) from the id\`{e}le group \(J_K\)
to the fractional ideal group \(I_K\) defined by, for each \(\alpha=(\alpha_P)\in J_K\),
\begin{equation}
      \pi_f(\alpha) = \prod_{P\nmid\infty}P^{v_P(\alpha_P)} \in I_K.
\label{eqn:pi_f_alpha}
\end{equation}
The kernel of the homomorphism \(\pi_f: J_K\to I_K\)
is given by
\[
      J_K^{S_\infty} = \prod_{P\mid\infty}K_P^*\times\prod_{P\nmid\infty}U_P.
\]

For a given nonempty finite subset \(S_c\subseteq S_f\), let \(I_K^{S_c}\) be the subgroup of \(I_K\)
which is finitely generated by the prime ideals in \(S_c\).
If \(S_c = \emptyset\), then \(I_K^{S_c}\) is the trivial group,
which has only one element being \(O_K\).
Now, the composite of \(\pi_f: J_K\to I_K\)  and the canonical projection 
\(\pi_{S_c}: I_K\to I_K/I_K^{S_c}\) can define a group homomorphism
from \(J_K\) to \(I_K/I_K^{S_c}\) as follows
\begin{equation}
\pi_{f,S_c}: J_K\to I_K\to I_K/I_K^{S_c}.
\label{eqn:pi_f_S_c}
\end{equation}

Let \(S=S_\infty\cup S_c\). 
The kernel of the above surjective homomorphism \(\pi_{f,S_c}: J_K\to I_K/I_K^{S_c}\)
is given by
\begin{equation}
      J_K^S = \prod_{P\in S}K_P^*\times\prod_{P\notin S}U_P.
\label{eqn:J_K_S}
\end{equation}

By applying the first group isomorphism theorem to the homomorphism \(\pi_{f,S_c}\),
we obtain
\begin{equation}
       J_K/J_K^S 
   =  J_K/\Ker\pi_{f,S_c}
   \iso \Image\pi_{f,S_c}
   =  I_K/I_K^{S_c}.
\label{eqn:I_K_J_K_first}
\end{equation}

By the unique prime ideal factorization theorem \cite[Theorem~1.3.3]{bib:neukirch99a}
and its relation with exponential valuation \cite[p.~69]{bib:neukirch99a},
the group homomorphism \(f:K^*\to I_K\) sending \(a\in K^*\) to \(a O_K\in I_K\)
is the same as the composite of the diagonal embedding \(K^* \hookrightarrow J_K\)
and the group homomorphism \(\pi_f: J_K\to I_K\) given in (\ref{eqn:pi_f_alpha}).
Therefore, the image of the following composite homomorphism
\begin{equation}
\hat{f}: K^* \hookrightarrow J_K\to I_K\to I_K/I_K^{S_c}
\label{eqn:phi_one}
\end{equation}
is the same as the image of the composite homomorphism 
of \(f:K^*\to I_K\) and \(\pi_{S_c}: I_K\to I_K/I_K^{S_c}\),
namely, \(K^* \to I_K\to I_K/I_K^{S_c}\).
Hence, the image of \(\hat{f}\)
is identical to the image of the following composite homomorphism
\begin{equation}
\tilde{f}: P_K \hookrightarrow I_K\to I_K/I_K^{S_c}
\label{eqn:phi_two}
\end{equation}
where \(P_K\) is the image of \(f:K^*\to I_K\),
and \(P_K \hookrightarrow I_K\) is the canonical injection.
That is to say \(\Image \hat{f} =\Image \tilde{f}\).

Noting that the homomorphism \(\hat{f}: K^* \to I_K/I_K^{S_c}\) given in (\ref{eqn:phi_one})
is also a composite of the diagonal embedding \(K^* \hookrightarrow J_K\)
and the homomorphism \(\pi_{f,S_c}: J_K \to I_K/I_K^{S_c}\) 
given in~(\ref{eqn:pi_f_S_c}),
we have \(\Ker\hat{f} = K^*\cap \Ker \pi_{f,S_c} = J_K^S  \cap K^*\).
It is clear that \(\Ker \tilde{f} = P_K\cap \Ker \pi_{S_c} = I_K^{S_c} \cap P_K\).
Then, by using the first group isomorphism theorem, we get
\begin{equation}
         K^*/(J_K^S \cap K^*) 
    =   K^*/\Ker\hat{f}
    \iso  \Image \hat{f} 
    =     \Image \tilde{f}
    \iso  P_K/\Ker \tilde{f}
    =  P_K/(I_K^{S_c}\cap P_K).
\label{eqn:I_K_J_K_second}
\end{equation}

We define \(C_K^S = J_K^S K^*/K^*\). 
Since \(K^*\subseteq J_K^S K^*\subseteq J_K\), 
we know that \(C_K^S\) is a normal subgroup of \(J_K/K^* = C_K\).
Similarly, we let \(Cl_K^{S_c} = I_K^{S_c} P_K/P_K\)
and, by \(P_K\subseteq I_K^{S_c} P_K\subseteq I_K\), 
we see that \(Cl_K^{S_c}\) is a normal subgroup of \(I_K/P_K = Cl_K\).
If \(S_c = \emptyset\), then \(Cl_K^{S_c}\) is the trivial group.
There is a group isomorphism as shown below.

\begin{PP}
For any finite subset \(S_c\subseteq S_f\), let \(S=S_\infty\cup S_c\). Then, we have
\begin{equation}
     C_K/C_K^S \iso Cl_K/Cl_K^{S_c}.
\label{eqn:Cl_K_C_K_iso}
\end{equation}
\label{prop:Cl_K_C_K}
\end{PP}

{\em Proof:}
By the second and third group isomorphism theorems, we obtain
\(
        C_K/C_K^S 
      = (J_K/K^*)/(J_K^S K^*/K^*)
   \iso J_K/J_K^S K^*
   \iso (J_K/J_K^S)/(J_K^S K^*/J_K^S)
   \iso (J_K/J_K^S)/(K^*/(J_K^S\cap K^*))
   \iso (I_K/I_K^{S_c})/(P_K/(I_K^{S_c}\cap P_K))   
   \iso (I_K/I_K^{S_c})/(I_K^{S_c} P_K/I_K^{S_c})      
   \iso I_K/I_K^{S_c} P_K
   \iso Cl_K/Cl_K^{S_c},
\)
where we have used (\ref{eqn:I_K_J_K_first}) and (\ref{eqn:I_K_J_K_second}).
The proof is finished. 
\eof

It follows from the finiteness of the ideal class group \(Cl_K\) 
and the group isomorphism given in~(\ref{eqn:Cl_K_C_K_iso}) 
that \(C_K^S\) is a subgroup of \(C_K\) of {\em finite} index.
According to \cite[Proposition~6.1.8]{bib:neukirch99a},
the subgroup \(C_K^S\) of \(C_K\) satisfies the topological requirement 
in the assumption of the existence theorem in the global class field theory
(see, \cite[Theorem~6.6.1]{bib:neukirch99a} and \cite[p.~172]{bib:tate67a}).
Therefore, there exists a unique finite abelian extension of~\(K\), denoted by \(H_{K,S_c}\),
corresponding to the subgroup \(C_K^S\) of \(C_K\).
The Galois group of the abelian extension \(H_{K,S_c}/K\) is given by
\[
    \Gal(H_{K,S_c}/K) \iso C_K/C_K^S \iso Cl_K/Cl_K^{S_c}.
\]
The number field \(H_{K,S_c}\) is called the {\em \(S_c\)-Hilbert class field} of \(K\).
If \(S_c = \emptyset\), then the subgroup \(C_K^S = C_K^{S_\infty} = J_K^{S_\infty} K^*/K^*\)
and \(H_{K,S_c} = H_K\), the Hilbert class field of \(K\).

For a finite place \(P\nmid\infty\),
the multiplicative group \(K_P^*\) can be embedded into \(J_K\) 
by the local embedding \(K_P^* \hookrightarrow J_K\)
sending \(a\in K_P^*\) to \(\alpha=(\alpha_P)\in J_K\) with 
\(\alpha_P = a\) and \(\alpha_{P'} = 1\) for all \(P'\in\Pfld_K\) with \(P'\neq P\).
Under the local embedding \(K_P^* \hookrightarrow J_K\),
it follows from (\ref{eqn:J_K_S}) and \(U_P\subseteq K_P^*\)
that \(K_P^*\subseteq J_K^S K^*\) for all \(P\in S_c\subseteq S\)
and \(U_P\subseteq J_K^S K^*\) for all \(P\nmid\infty\)
for the subgroup \(C_K^S = J_K^S K^*/K^*\) of \(C_K\).

By using the criteria given in \cite[Proposition~2.5.3]{bib:niederreiter01a}, 
which apply for all {\em non-archimedean} places of global fields
including {\em finite} places of algebraic number fields,
we see that the \(S_c\)-Hilbert class field \(H_{K,S_c}\) of \(K\)
has the following two properties

1) {\em Each finite place  \(P\nmid\infty\) of \(K\) is unramified in \(H_{K,S_c}\)},
and 

2) {\em Each finite place \(P\in S_c\) splits completely in \(H_{K,S_c}\)}.

It is ready to build a tower of fields as follows
\[
      K=K_1^{(S_1)} \subseteq K_2^{(S_2)}\subseteq K_3^{(S_3)}\subseteq\cdots
\]
where \(S_1 = S_c\) and for \(i\in\nfld\), \(K_{i+1}^{(S_{i+1})} = H_{K_i^{(S_i)},S_i}\), 
the \(S_i\)-Hilbert class field of \(K_i^{(S_i)}\), 
and \(S_{i+1}\) is the set of all finite places of \(K_{i+1}^{(S_{i+1})}\) 
lying over the places in \(S_i\).
This tower of fields is called the {\em \(S_c\)-Hilbert class field tower} of \(K\).
If \([K_i^{(S_i)}:K]\to\infty\) as \(i\to\infty\),
we say that \(K\) has an {\em infinite} \(S_c\)-Hilbert class field tower.

Let \(K^S = K^*\cap J_K^S = \{a\in K^*\mid v_P(a) = 0 \mbox{ for all } P\in S_f \setminus S_c\}\), 
each element of which is called an {\em \(S\)-unit} of \(K\).
If \(S_c = \emptyset\), i.e., \(S=S_\infty\), 
then \(K^S = K^{S_\infty}= K^*\cap J_K^{S_\infty} = U_K\),
the group of units of \(O_K\).
The generalized Golod--Shafarevich condition has been given by Ihara \cite[p.~705]{bib:ihara83a},
which asserts that if, for some prime \(p\),
\begin{equation}
d_p (Cl_K/Cl_K^{S_c}) \geq 2 + 2 \sqrt{d_p K^S + 1}
\label{eqn:GS_condition_generalized}
\end{equation}
then \(K\) has an infinite \(S_c\)-Hilbert class field tower.
If \(S_c = \emptyset\), then \(Cl_K/Cl_K^{S_c}\iso Cl_K\) and \(K^S = U_K\) 
and hence that
(\ref{eqn:GS_condition_generalized}) reduces to (\ref{eqn:GS_condition_Cl_U}).

It was mentioned in \cite[p.~705]{bib:ihara83a} that 
the proof presented in Roquette \cite{bib:roquette67a} for the case \(S_c = \emptyset\)
can be generalized immediately
to the nonempty case \(S_c \neq \emptyset\) 
just by considering the elements of \(S_c\) as additional ``infinite primes''.
A proof for an even more general theorem including the generalized Golod--Shafarevich condition
(\ref{eqn:GS_condition_generalized}) has been briefly presented by Maire \cite{bib:maire96a}.
The generalized Golod--Shafarevich condition (\ref{eqn:GS_condition_generalized})
for an infinite class field tower of an algebraic number field
is an analogue of the one applied for an infinite class field tower
of an algebraic function field over a finite field
(see Schoof \cite[p.~11]{bib:schoof90a}, Niederreiter and Xing \cite[p.~59]{bib:niederreiter01a},
and the monograph by Serre \cite{bib:serre2020c}).

For the purpose of applicability, we are going to present a slight generalization of
the Golod--Shafarevich condition (\ref{eqn:GS_condition_Cl_U}),
which will be suitable for application to the code constructions in this paper.

Let \(L/K\) be a Galois extension of algebraic number fields.
We denote by \(P(L/K)\) the set of all nonzero prime ideals of \(K\) which split completely in \(L\).
According to Chebotarev's density theorem, it is known that the Dirichlet density of the set \(P(L/K)\)
is equal to~\(1/[L:K]\), where \([L:K ]\) is the degree of extension (see \cite[p.~547]{bib:neukirch99a}).

If \(L=H_K\), the Hilbert class field of \(K\), then the set  \(P(L/K)\) consists precisely of
the nonzero prime ideals of \(K\) which are also principal ideals.
Since \([H_K:K]=h_K\), the class number of \(K\), the Dirichlet density of the set \(P(H_K/K)\) is equal to
\(1/h_K\). By the finiteness of \(h_K\), the set \(P(H_K/K)\) must have positive density.
Therefore, it is {\em guaranteed} that one can choose an {\em arbitrarily large}
nonempty finite subset \(S_c \subseteq P(H_K/K)\).

For any nonempty finite subset \(S_c \subseteq P(H_K/K)\),
each of the prime ideals in \(S_c\) has the form of \(aO_K\) with \(a\in O_K\) and \(a\neq 0\),
which is a principal ideal of \(O_K\).
Hence, \(S_c\subseteq P_K\), the subgroup of nonzero principal fractional ideals of~\(I_K\).
Noting that \(I_K^{S_c}\) is the subgroup of \(I_K\)
that is finitely generated by the prime ideals in~\(S_c\),
we have \(I_K^{S_c}\subseteq P_K\).
Thus, \(I_K^{S_c} P_K = P_K\). 
Consequently, \(Cl_K^{S_c} = I_K^{S_c} P_K/P_K\) is the trivial group.

Therefore, for any finite subset \(S_c \subseteq P(H_K/K)\), it is true that \(Cl_K/Cl_K^{S_c}\iso Cl_K\)
and then the \(p\)-rank \(d_p(Cl_K/Cl_K^{S_c}) = d_p Cl_K\) for any prime \(p\).
By using (\ref{eqn:GS_condition_generalized}), we have the following criterion.

\begin{PP}
Let \(S_c \subseteq P(H_K/K)\) be a finite subset and \(S=S_\infty\cup S_c\). 
If, for some prime \(p\),
\begin{equation}
d_p (Cl_K) \geq 2 + 2 \sqrt{d_p K^S + 1}
\label{eqn:p_S_Hilbert}
\end{equation}
then \(K\) has an infinite \(S_c\)-Hilbert class field tower.
\label{prop:p_S_Hilbert}
\end{PP}

The \(p\)-rank of \(K^S\), i.e., \(d_p K^S\), in (\ref{eqn:p_S_Hilbert}) can be calculated as follows.
According to a generalized Dirichlet's unit theorem \cite[Proposition~6.1.1]{bib:neukirch99a}, 
we know that \(K^S\) is a finitely generated abelian group and that
\begin{equation}
    d_p K^S = |S| - 1 + \delta_p(K) = |S_c| + |S_\infty| - 1 + \delta_p(K)
\label{eqn:Dirichlet_unit_generalized}
\end{equation}
where \(\delta_p(K)\) equals \(1\) if \(K\) contains a primitive \(p\)th root of unity 
and \(0\) otherwise.

By employing an infinite tower of \(S_c\)-Hilbert class fields
with a constraint on the norms of prime ideals in \(S_c\), 
we can obtain a lower bound on the quantity \(A(r,q)\).
Subsequently, by using Lenstra codes,
we can achieve a lower bound on the rate function~\(R_q(\delta)\), as shown below.

\begin{PP}
Let \(K\) be a number field of discriminant \(\Delta_K\) with \(|\Delta_K| > 1\).
Let \(r\) and~\(q\) be two positive integers satisfying \(2\leq r\leq q\) and
\(S_c \subseteq P(H_K/K)\) be a nonempty finite subset, each prime ideal in which has norm between \(r\) and \(q\).
If, for some prime \(p\),
\begin{equation}
d_p (Cl_K) \geq 2 + 2 \sqrt{|S_c| + |S_\infty| + 1}
\label{eqn:construction_criterion}
\end{equation}
then
\begin{equation}
   R_q(\delta) \geq (1-\delta)\frac{\log r}{\log q} - \frac{\log\sqrt{|\Delta_K|}}{|S_c|\log q}
\label{eqn:Rq_CFT_bound}
\end{equation}
for \(0 < \delta < 1\).
\label{prop:connection_CFT_Lenstra}
\end{PP}

{\em Proof:} 
It is seen from (\ref{eqn:Dirichlet_unit_generalized}) that
\(d_p K^S \leq |S_c| + |S_\infty|\).
If, for some prime \(p\), the condition (\ref{eqn:construction_criterion}) is satisfied,
then (\ref{eqn:p_S_Hilbert}) holds.
According to Proposition~\ref{prop:p_S_Hilbert}, \(K\) has an infinite \(S_c\)-Hilbert class field tower.

Let \(K=K_1^{(S_1)} \subseteq K_2^{(S_2)}\subseteq\cdots\) be 
the infinite \(S_c\)-Hilbert class field tower of~\(K\), where \(S_1 = S_c\).
Then, \(\lim_{i\to\infty}[K_i^{(S_i)}:K]\to\infty\). 
Moreover, \(K_{i+1}^{(S_{i+1})}\) is the \(S_i\)-Hilbert class field of \(K_i^{(S_i)}\) 
for \(i\in\nfld\).

Let \(K_i = K_i^{(S_i)}\) and \(\Delta_i = |\Delta_{K_i}|\) for \(i\in\nfld\).
Then,  \([K_i:K]\to\infty\) as \(i\to\infty\).
Since each finite place \(P\nmid\infty\) of \(K_i\) is unramified in \(K_{i+1}\),
by virtue of an analogue of the Riemann--Hurwitz formula for algebraic number fields
\cite[Proposition 3.3.13]{bib:neukirch99a},
we have 
\[
    \log\sqrt{\Delta_{i+1}} = [K_{i+1}:K_i] \log\sqrt{\Delta_i}
\]
for \(i\in\nfld\). Hence, 
\begin{equation}
    \log\sqrt{\Delta_i} = [K_i:K_1] \log\sqrt{|\Delta_K|}
\label{eqn:log_Delta_K}
\end{equation}
for \(i\in\nfld\), which implies \(\Delta_i\to\infty\) as \(i\to\infty\).

For each nonzero prime ideal \(P_i\in S_i\), it {\em splits completely} in \(K_{i+1}\)
for \(i\in\nfld\).
Therefore, there are \([K_{i+1}:K_i]\) nonzero prime ideals in \(S_{i+1}\)
lying over each \(P_i\in S_i\) for \(i\in\nfld\).
Consequently, we have
\(
      |S_{i+1}| = [K_{i+1}:K_i]|S_i|
\)
for \(i\in\nfld\). 
Moreover, for \(i\in\nfld\) with \(i\geq 2\), 
each \(P_i\in S_i\) has the same norm as 
that of the prime ideal \(P\in S_c\) lying below~\(P_i\).
Hence, for \(i\in\nfld\),
the norm of each prime ideal \(P_i\in S_i\) is between \(r\) and \(q\).
Then, it follows that
\begin{equation}
    N_{r,q}(K_i) \geq |S_i| = [K_i:K_1]|S_1| = [K_i:K_1]|S_c|
\label{eqn:N_r_q_S_c}
\end{equation}
for \(i\in\nfld\).

Therefore, we get
\begin{equation}
     A(r,q)
  = \limsup_{\Delta\rightarrow\infty}\frac{N_{r,q}(\Delta)}{\log\sqrt{\Delta}}
  \geq \limsup_{i\rightarrow\infty}\frac{N_{r,q}(\Delta_i)}{\log\sqrt{\Delta_i}}
  \geq \limsup_{i\rightarrow\infty}\frac{N_{r,q}(K_i)}{\log\sqrt{\Delta_i}}
  \geq \frac{|S_c|}{\log\sqrt{|\Delta_K|}}
\label{eqn:A_r_q_CFT_bound}
\end{equation}
where the last inequality follows from (\ref{eqn:log_Delta_K}) and (\ref{eqn:N_r_q_S_c}).
It is clear that (\ref{eqn:Rq_CFT_bound}) can be deduced from
(\ref{eqn:Rq_A_r_q_lower_bound}) in Proposition~\ref{prop:Rq_A_r_q}
and (\ref{eqn:A_r_q_CFT_bound}).

The proof of this proposition is complete.
\eof

Proposition~\ref{prop:connection_CFT_Lenstra} is the theoretical foundation 
for the following practical code constructions in this paper.
For the application, a lower bound for the \(p\)-rank of the ideal class group \(Cl_K\), i.e., \(d_p Cl_K\),
in the Golod--Shafarevich condition (\ref{eqn:construction_criterion})
will be needed.


\section{An applicable lower bound for the rate function}
\label{sec:lowerbound}

Based on theory of Lenstra codes from geometry of numbers given in Section~\ref{sec:Lenstra}
and theory of infinite towers of Hilbert class fields in Section~\ref{sec:Hilbert},
we can obtain an applicable lower bound on the rate function~\(R_q(\delta)\) in the following.

Let \(p_i\in\Pfld\) be the \(i\)th prime for \(i\in\nfld\).
For an odd prime \(p\) and an integer \(a\in\Zfld\),
the Legendre symbol \(\big(\frac{a}{p}\big)\in\{-1,0,1\}\)
is function of \(a\) and \(p\) defined as \(0\) when \(p\mid a\) and \(\pm 1\) 
depending on whether \(a\) is a quadratic residue modulo \(p\) when \(p\nmid a\).

\begin{TT}
Let \(q\in\nfld\) with \(q\geq 2\).
If there exist \(r,\ell,k\in\nfld\) which satisfy the following three conditions

1) \(2\leq r\leq q\),

2) \(k+2\leq \frac{1}{4}(\ell-2)^2-(\ell-2)\), and

3) There are at least \(2k\) primes \(p\) 
   satisfying \(r\leq p^2\leq q\),  \(p>p_\ell\), and 
   \(p \equiv 3 \pmod{4}\),

\noindent then 
\begin{equation}
   R_q(\delta) \geq (1-\delta)\frac{\log r}{\log q} - \frac{\log\sqrt{D}}{k\log q}
\label{eqn:Rq_general_bound}
\end{equation}
for \(0 < \delta < 1\), where \(D = 4 p_1\cdots p_\ell\).
\label{thrm:liang_applicable}
\end{TT}

{\em Proof:} 
We consider the quadratic number field \(K = \Qfld(\sqrt{\Delta_K})\),
where \(\Delta_K = D\) or \(-D\). 
It is clear that the discriminant of \(K\) is \(\Delta_K\).
If \(\Delta_K > 0\), then \(K\) is a real quadratic field.
If \(\Delta_K < 0\), then \(K\) is an imaginary quadratic field.

Since \(\Delta_K\) is the discriminant of the quadratic field \(K\),
it is a {\em fundamental discriminant}.
Thus, according to \cite[Theorem~2.2.15]{bib:cohen07a},
the Kronecker symbol \(\big(\frac{\Delta_K}{n}\big)\), \(n\in \Zfld\), 
defines a real primitive Dirichlet character modulo \(|\Delta_K|=D\).
We note that if \( n = p\), \(p\) an odd prime, then the Kronecker symbol \(\big(\frac{\Delta_K}{n}\big)\)
has the same value as that of the Legendre symbol \(\big(\frac{\Delta_K}{p}\big)\).
Therefore, if \(p\) is an odd prime, then
\[
    \left(\frac{-D}{p}\right)
 =
    \left(\frac{-1}{p}\right) \left(\frac{D}{p}\right).
\]
It is known, by Euler’s criterion, that \(\big(\frac{-1}{p}\big)=-1\) for every prime \(p \equiv 3 \pmod{4}\).
Thus, if  \(p \equiv 3 \pmod{4}\) and \(p\nmid D\),
then either \(\big(\frac{D}{p}\big)=-1\) or \(\big(\frac{-D}{p}\big)=-1\).

For each of the primes \(p \equiv 3 \pmod{4}\) in condition~3), we have \(p>p_\ell\)
and thus \(p\nmid D\). Therefore, each of the primes \(p\) in condition~3)
satisfies either \(\big(\frac{D}{p}\big)=-1\) or \(\big(\frac{-D}{p}\big)=-1\).
Based on this, condition~3) can be reduced to the following two cases, namely, either

{\em Case 1)}  There are at least \(k\) primes \(p \equiv 3 \pmod{4}\)
   satisfying \(r\leq p^2\leq q\),  \(p>p_\ell\), and  \(\big(\frac{D}{p}\big)=-1\), or

{\em Case 2)}  There are at least \(k\) primes \(p \equiv 3 \pmod{4}\)
   satisfying \(r\leq p^2\leq q\),  \(p>p_\ell\), and  \(\big(\frac{-D}{p}\big)=-1\).

Now, we can choose the quadratic number field \(K = \Qfld(\sqrt{\Delta_K})\),
where \(\Delta_K = D\) in Case~1) and  \(\Delta_K = -D\) in Case~2).
It follows from condition~3) that
there are at least \(k\) primes \(p \equiv 3 \pmod{4}\)
satisfying \(r\leq p^2\leq q\),  \(p>p_\ell\), and  \(\big(\frac{\Delta_K}{p}\big)=-1\).
By using \cite[Proposition~3.4.3]{bib:cohen07a}, we know that
each of these \(k\) primes \(p\) is {\em inert} in \(K\),
namely, the principal ideal \(pO_K\) being a nonzero prime ideal of \(O_K\).
Thus, the set
\(
   \{ pO_K\mid r\leq p^2\leq q,\; p>p_\ell,\mbox{ and } \big(\frac{\Delta_K}{p}\big) = -1\}
\)
has cardinality of at least \(k\), and it is a subset of \(P(H_K/K)\).

Therefore, condition~3) implies that there exists a subset, denoted by \(S_c\),
of the set
\(
   \{ pO_K\mid r\leq p^2\leq q,\; p>p_\ell,\mbox{ and } \big(\frac{\Delta_K}{p}\big) = -1\}
\)
such that \(|S_c| = k\) and \(S_c \subseteq P(H_K/K)\).

Note that the norm of each \(pO_K\in S_c\) is given by
\[
   \norm(pO_K) = |\NKQ(p)| = p^{[K:\Qfld]} = p^2.
\]
Thus, every nonzero prime ideal in \(S_c\) has norm between \(r\) and \(q\).

In the following, we are going to verify that \(K\) and \(S_c\) satisfy the generalized Golod--Shafarevich condition (\ref{eqn:construction_criterion}) for \(p = 2\).
According to Proposition~\ref{prop:connection_CFT_Lenstra}, 
(\ref{eqn:Rq_general_bound}) will then be deduced immediately from (\ref{eqn:Rq_CFT_bound}).

For any given prime \(p\), the \(p\)-rank \(d_p Cl_K\) has a lower bound as follows
(see Schoof \cite{bib:schoof86a} and Martinet \cite{bib:martinet78a}).
Let \(K/\Qfld\) be a Galois extension and~\(\theta\) denote the number of finite places 
of \(\Qfld\) (i.e., primes of \(\Zfld\)) ramified in \(K/\Qfld\).
Assume that \(\Gal(K/\Qfld)\) is cyclic of order \(p\).
Then, we have
\(d_p Cl_K \geq \theta - 1 - \delta_p(\Qfld)\),
where \(\delta_p(\Qfld)\) equals \(1\) if \(\Qfld\) contains a primitive \(p\)th root of unity 
and \(0\) otherwise.

Taking \(p=2\), we obtain \(d_2 Cl_K \geq \theta - 2\) for the given quadratic field 
\(K = \Qfld(\sqrt{\Delta_K})\).
It is known that a finite place \(p\) in \(\Qfld\) is ramified in \(K/\Qfld\)
if and only if \(p\mid \Delta_K\) (see \cite[Proposition~3.4.3]{bib:cohen07a}).
Thus, \(\theta\) is equal to the number of primes 
which divide the discriminant \(\Delta_K\) of \(K\), i.e., \(\theta = \ell\).
Hence,
\begin{equation}
    d_2(Cl_K) \geq \ell - 2.
\label{eqn:d_2_Cl_K_S_c}
\end{equation}

The number of infinite places of \(K = \Qfld(\sqrt{\Delta_K})\) is given by
\(|S_{\infty}| = 2\) if \(\Delta_K > 0\) and by \(|S_{\infty}| = 1\) if \(\Delta_K < 0\).
Thus, 
\begin{equation}
    |S_c| + |S_\infty| \leq k + 2.
\label{eqn:d_2_K_S}
\end{equation}

By (\ref{eqn:d_2_Cl_K_S_c}),  (\ref{eqn:d_2_K_S}), and the given condition~2),
we see that \(K\) and \(S_c\) satisfy the generalized Golod--Shafarevich condition (\ref{eqn:construction_criterion}) for \(p=2\) as follows
\[
   d_2 (Cl_K) \geq 2 + 2 \sqrt{ |S_c| + |S_\infty| + 1}.
\]
Hence,  
(\ref{eqn:Rq_general_bound}) follows from (\ref{eqn:Rq_CFT_bound})
in Proposition~\ref{prop:connection_CFT_Lenstra}.
The proof is complete.
\eof

It is noticed that (\ref{eqn:Rq_general_bound}) in Theorem~\ref{thrm:liang_applicable}
is inferred from (\ref{eqn:Rq_CFT_bound}) in Proposition~\ref{prop:connection_CFT_Lenstra},
which in~turn is from (\ref{eqn:Rq_A_r_q_lower_bound}) in Proposition~\ref{prop:Rq_A_r_q}.

In the following, we consider how to employ Theorem~\ref{thrm:liang_applicable}
for \(q\in\nfld\) sufficiently large.

The lower bound for the rate function \(R_q(\delta)\) given in Proposition~\ref{prop:Rq_A_r_q},
namely, the right-hand side of (\ref{eqn:Rq_A_r_q_lower_bound}),
is related to a positive integer \(r\) with \(2\leq r\leq q\).
In order for this lower bound to be at least comparable to 
the Gilbert--Varshamov bound \(\RGV\) as given in~(\ref{eqn:RGV_asymptotic}),
one cannot choose \(r\) as small as \(O(q^{1-\ve})\) for some constant \(\ve\in(0,1)\).

Thus, we may choose \(r = \lceil\lambda q \rceil\) with an arbitrary but fixed number
\(\lambda\in(0,1]\).
If \(\lambda = 1\) is chosen, i.e., \(r = q\), then, from Lemma~\ref{lemma:A_q_q_bigO},
\(A(q,q) \leq 1/\hat{c}\) for a constant \(\hat{c}>0\). 
In this case, the right-hand side of (\ref{eqn:Rq_A_r_q_lower_bound}) is 
not greater than \(1 - \delta -\frac{\hat{c}}{\log q}\).
If we choose \(r = \lceil\lambda q \rceil\) with \(0<\lambda < 1\), the right-hand side of (\ref{eqn:Rq_A_r_q_lower_bound}) cannot exceed
\(1 - \delta -\frac{\tilde{c}}{\log q}\) with \(\tilde{c} = -(1-\delta)\log\lambda > 0\)
for \(0<\delta<1\). 
Therefore, with a choice of \(r = \lceil\lambda q \rceil\) for a constant \(\lambda\in(0,1]\),
the best attainable lower bound from~(\ref{eqn:Rq_A_r_q_lower_bound}) 
in Proposition~\ref{prop:Rq_A_r_q}
has asymptotically a {\em logarithmic} order in \(q\)
as the Gilbert--Varshamov bound in (\ref{eqn:RGV_asymptotic}).

As indicated by (\ref{eqn:hat_N_r_q}) defining the number \(\hat{N}_{r,q}(K)\),
in \cite{bib:lenstra86a} Lenstra employed an extended code from geometry of numbers, 
which takes infinite places into account,
but it was chosen in \cite{bib:lenstra86a} that
\(r = \lceil\frac{1}{2} q \rceil\) (i.e., \(\lambda = \frac{1}{2}\))
for nonzero prime ideals (i.e., finite places) of norm or its power between \(r\) and \(q\).
By using the previous analysis in the same spirit, 
we can conclude that 
an achievable lower bound obtained by Lenstra in \cite{bib:lenstra86a} 
for the rate function~\(R_q(\delta)\)
at~best has asymptotically a {\em logarithmic} order in \(q\) 
as does the Gilbert--Varshamov bound in~(\ref{eqn:RGV_asymptotic}),
even though the generalized Riemann hypothesis (GRH) was assumed in \cite{bib:lenstra86a}.

To beat the \(q\)-ary Gilbert--Varshamov bound \(\RGV\) 
by a better lower bound for large \(q\),
we find that {\em it is actually feasible} with the help of Theorem~\ref{thrm:liang_applicable}.
We can make use of (\ref{eqn:Rq_general_bound}) in Theorem~\ref{thrm:liang_applicable}
by choosing a sequence of numbers \(\lambda_q\), 
satisfying \(0<\lambda_q < 1\) with \(\lambda_q\to 1\) as \(q\to\infty\),
to define \(r=r(q)=\lceil\lambda_q q \rceil\),
while dealing with carefully the trade-off between choice of \(r(q)\) (i.e., \(\lambda_q\))
and that of \(\ell = \ell(q)\) and \(k=k(q)\)
to maximize the attainable lower bound from (\ref{eqn:Rq_general_bound}).
This approach will need a subtle analysis of \(k=k(q)\),
given \(r=r(q)\) and \(\ell = \ell(q)\),
in condition~3) of Theorem~\ref{thrm:liang_applicable},
and it relies on analytic number theory
for the distribution of prime numbers.

In contrast to previous code constructions which have been based on
{\em a single} infinite tower of number fields in \cite{bib:lenstra86a}
or of algebraic function fields over finite fields in \cite{bib:tsfasman82a}, 
\cite{bib:serre2020c}, \cite{bib:schoof90a}, 
\cite{bib:garcia95a}, and \cite{bib:bassa14a},
the code construction in the proof of Theorem~\ref{thrm:liang_applicable} uses
 {\em two} differentiated but complementary towers of Hilbert class fields
 with at least one being infinite.

It is ready now to establish proofs of 
Theorem~\ref{thrm:thrm_liang_constant_C} and Theorem~\ref{thrm:liang},
as will be presented in the subsequent sections. 
It is emphasized that the GRH is {\em not} assumed in the proofs,
and each new lower bound in Theorem~\ref{thrm:thrm_liang_constant_C} and Theorem~\ref{thrm:liang}
is an {\em unconditional} result.


\section{Proof of Theorem~\ref{thrm:thrm_liang_constant_C} }
\label{sec:ProofLiang1over6}

We apply Theorem~\ref{thrm:liang_applicable} for \(q\in\nfld\) sufficiently large,
which needs an appropriate choice of the three positive integers \(r,\ell,k\in\nfld\), depending on \(q\).

First, we choose
\begin{equation}
    \ve_q = \frac{(\log q)\log\log q}{C_0 q^{1/6}}
\label{eqn:liang_1over6_ve_q_RGH}
\end{equation}
where \(C_0 > 1\) is a suitable absolute constant such that 
\(0<\ve_q<1\) for all \(q\geq 27\).
It is seen that \(\ve_q\to 0\) as \(q\to\infty\).
We define
\begin{equation}
     r = r(q) = \lceil(1 - \ve_q)^2 q \rceil
\label{eqn:liang_1over6_r_q_ceil_GRH}
\end{equation}
for \(q\geq 27\).

Next, we set
\[
  x = \sqrt{q},  \; \; \;  y = (1-\ve_q) \sqrt{q}
\]
and choose
\begin{equation}
     \ell = \ell(q) = \lfloor x^{1/3}\rfloor
\label{eqn:liang_1over6_ell_q_x_log_x_GRH}
\end{equation}
for \(q\geq 27\).

Finally, we define 
\begin{equation}
    k = k(q) = \left\lfloor\tfrac{1}{4}(\ell-2)^2-(\ell-2)\right\rfloor - 2
\label{eqn:liang_1over6_k_q_ell_q_GRH}
\end{equation}
for \(q\in\nfld\) such that \(\ell(q)\geq 10\).

The three positive integers \(r,\ell,k\in\nfld\), depending on \(q\), 
as chosen in the above satisfy condition~1) and condition~2)  of Theorem~\ref{thrm:liang_applicable}
for \(q\) large enough. 
In the following, we need to check that
condition~3) of Theorem~\ref{thrm:liang_applicable} is also satisfied by the chosen \(r,\ell,k\in\nfld\)
provided that \(q\) is sufficiently large. 

Let \(\pi(x;4,3)\) denote the number of primes less than or equal to \(x\)
which are congruent to \(3 \bmod 4\).

We are going to prove that the following cardinality of a subset
\[
   N_q
\dfnt
 \left|\{ p\in\Pfld \mid \sqrt{r} \leq p \leq \sqrt{q}, \;  p>p_\ell,\mbox{ and } p \equiv 3 \!\!\!\! \pmod {4}\}\right|
\]
satisfies \(N_q \geq 2k\) for \(q\) large enough. 

Noting that \(\ell = \lfloor x^{1/3}\rfloor = \lfloor q^{1/6} \rfloor\), we get
\[
      \sqrt{r} 
   \geq (1-\ve_q) \sqrt{q}
   = (1-\ve_q) x
   > 0.99 x
   > 2 x^{1/3}\log(x^{1/3})
   \geq 2 \ell\log\ell
   > p_\ell
\]
for \(q\) large enough.

Therefore, we have
\[
   N_q
 = \left|\{ p\in\Pfld \mid \sqrt{r} \leq p \leq \sqrt{q},\; p \equiv 3 \!\!\!\! \pmod {4}\}\right|
 = \left|\{ p\in\Pfld \mid y \leq p \leq x,\; p \equiv 3 \!\!\!\! \pmod {4}\}\right|.
\]
Hence, for \(q\) sufficiently large,
\begin{equation}
N_q \geq \pi(x;4,3)-\pi(y;4,3)
\label{eqn:liang_1over6_N_q_primes}
\end{equation}
where \(x = \sqrt{q}\) and \(y = (1-\ve_q) x\).

For the given functions \(f(q)\) and \(g(q)\) of \(q\in\nfld\),
if \(\lim_{q\to\infty} \frac{f(q)}{g(q)} = 1\),
we say that the two functions \(f(q)\) and \(g(q)\) 
are asymptotically equivalent, denoted by \(f(q)\sim g(q)\), as \(q\to\infty\).
By the prime number theorem, we know that the positive integer \(D = 4 p_1\cdots p_\ell\) 
defined in Theorem~\ref{thrm:liang_applicable} satisfies
\begin{equation}
     \log D = \log 4 + \sum_{p\leq p_\ell} \log p \sim p_\ell \sim \ell\log \ell
\label{eqn:liang_1over6_log_D_ell_log_ell}
\end{equation}
as \(q\to\infty\) (see \cite[Chapter~4]{bib:apostol76a}).
By (\ref{eqn:liang_1over6_ell_q_x_log_x_GRH}), we have 
\(\ell\sim q^{1/6}\)
and hence
\(\log\ell\sim\frac{1}{6}\log q\)
and
\(p_\ell\sim \ell\log \ell \sim\frac{1}{6}q^{1/6}\log q\) as \(q\to\infty\).

We define \(h=x-y\) with \(x = \sqrt{q}\) and \(y = (1-\ve_q) x\), which satisfies
\begin{equation}
   h = \ve_q x  \sim C_1 x^{2/3} (\log x)\log\log x
\label{eqn:liang_1over6_h_length}
\end{equation}
for \(q\) sufficiently large with a suitable constant \(C_1 > 0\).
 
By using theory on the distribution of primes in arithmetic progressions
and in short intervals
(see, e.g., \cite[Theorem~3]{bib:baker1997a}, \cite[p.~656, Theorem]{bib:motohashi71a}),
we know that, for any fixed real number \( C_2 \in (0,1) \),
\begin{equation}
    \pi(x;4,3)-\pi(y;4,3)
 =\pi(x;4,3)-\pi(x-h;4,3)
 \geq
   C_2\frac{h}{2\log x}
\label{eqn:liang_1over6_primes_in_short_interval}
\end{equation}
for \(x= \sqrt{q}\) sufficiently large.
Note that (\ref{eqn:liang_1over6_primes_in_short_interval}) holds for all \(h\) 
with \(x^{\beta}\leq h\leq x(\log x)^{-1}\) where \(\beta\in\rfld\) is any constant satisfying \(5/8 < \beta < 1\),
and it is a generalization of Ingham's theorem
on the distribution of primes in short intervals \cite{bib:ingham37a}.

It follows from (\ref{eqn:liang_1over6_N_q_primes}),
(\ref{eqn:liang_1over6_primes_in_short_interval}), and (\ref{eqn:liang_1over6_h_length}) that
\[
     N_q 
  \geq
    C_2\frac{h}{2\log x}
  \sim 
   C_3 x^{2/3}\log\log x
\]
for \(q\) sufficiently large with a suitable constant \(C_3 > 0\).
On the other hand, we have
\[
  2k \sim \frac{1}{2}\ell^2 \sim \frac{1}{2}x^{2/3}
\]
for  \(q\) sufficiently large. 
Hence,
\(N_q \geq 2k\) for \(q\) large enough. 

In summary, we have shown that
the three positive integers \(r,\ell,k\in\nfld\)
as chosen in (\ref{eqn:liang_1over6_r_q_ceil_GRH}),
(\ref{eqn:liang_1over6_ell_q_x_log_x_GRH}), and (\ref{eqn:liang_1over6_k_q_ell_q_GRH}),
depending on \(q\), satisfy all three conditions
in Theorem~\ref{thrm:liang_applicable} for \(q\in\nfld\) sufficiently large.

Now, by using (\ref{eqn:Rq_general_bound}), we obtain
\begin{equation}
        R_q(\delta) 
   \geq (1-\delta)\frac{\log r}{\log q} - \frac{\log\sqrt{D}}{k\log q}
   =  1 - \delta - (1-\delta)\frac{\log(q/r)}{\log q} - \frac{\log D}{2k\log q}
\label{eqn:liang_1over6_R_q_delta_apply_GRH}
\end{equation}
for \(q\geq C_4\) and \(0 < \delta < 1\), where \(C_4>0\) is an absolute constant.

It is clear from (\ref{eqn:liang_1over6_r_q_ceil_GRH}) that
\begin{equation}
         \log\frac{q}{r}
    \leq \log\frac{q}{(1 - \ve_q)^2 q}
       = 2 \log\frac{1}{1 - \ve_q}
       < 3 \ve_q
\label{eqn:liang_1over6_frac_q_r_GRH}
\end{equation}
for \(q\) sufficiently large.

Moreover, we can see that
\[
           \frac{\log D}{4k} 
      \sim \frac{\log D}{\ell^2}
      \sim \frac{\log \ell}{\ell}
      \sim \frac{\log x}{3x^{1/3}}
      \sim \frac{C_0\ve_q}{C_5\log\log q}
\]
as \(q\to\infty\), 
where \(C_5 > 0\) is some constant.

Hence, for \(q\) large enough,
\begin{equation}
           \frac{\log D}{2k} 
         < 2 C_0\ve_q.
\label{eqn:liang_1over6_frac_logD_k_GRH}
\end{equation}

It follows from (\ref{eqn:liang_1over6_R_q_delta_apply_GRH}), (\ref{eqn:liang_1over6_frac_q_r_GRH}),
and (\ref{eqn:liang_1over6_frac_logD_k_GRH}) that there exists an absolute constant \(C_6 > 0\) such that,
for all \(q\geq C_6\) and \(0<\delta<1\),
\[
      R_q(\delta) 
   >  1 - \delta - \frac{3\ve_q}{\log q} - \frac{2 C_0\ve_q}{\log q}
   >  1 - \delta - \frac{5 C_0\ve_q}{\log q}
   =  1 - \delta - 5 \frac{\log\log q}{q^{1/6}}.
\]

Consequently, there exists an absolute constant \(C > 5\) such that
\[
      R_q(\delta) 
   >  
     1 - \delta - C \frac{\log\log q}{q^{1/6}}
\]
for all \(q\geq 3\) and \(0<\delta<1\).

Therefore, the proof of Theorem~\ref{thrm:thrm_liang_constant_C}  is complete.
\eof

\begin{RR}
The estimate of (\ref{eqn:liang_1over6_primes_in_short_interval}) can be proved directly
as a straightforward generalization of 
Ingham's proof in \cite{bib:ingham37a}.
In this scenario, the Dirichlet \(L\)-function \(L(s,\chi_4)\) of \(s\in\cfld\) will be involved,
where \(\chi_4\) is the primitive Dirichlet character modulo 4.
It is noted that \(L(s,\chi_4)\) has no Landau--Siegel zero.
Then, the explict formula (see, e.g., \cite{bib:davenport80a})
and the zero-density estimates for \(L(s,\chi_4)\) 
(see, e.g., \cite{bib:montgomery71a})
will enable us to prove (\ref{eqn:liang_1over6_primes_in_short_interval}).
In this way, the absolute constant \(C_0>0\) for sufficiently large \(x\geq C_0\) satisfying (\ref{eqn:liang_1over6_primes_in_short_interval})
can be effectively computed.
\end{RR}


\section{Proof of Theorem~\ref{thrm:liang}}
\label{sec:ProofLiang}

In the following, we apply Theorem~\ref{thrm:liang_applicable} for all \(q\in\nfld\) with \(q \geq Q\dfnt \exp(29)\).

A suitable choice of the three positive integers \(r,\ell,k\in\nfld\), depending on \(q\), is given as follows.
First, we choose
\begin{equation}
    \ve_q = \frac{1}{(\log q)^{1/3}}.
\label{eqn:thrm_liang_ve_q}
\end{equation}
It is clear that \(0<\ve_q<1\) for all \(q\geq Q\) and \(\ve_q\to 0\) as \(q\to\infty\).
We define
\begin{equation}
     r = r(q) = \lceil(1 - \ve_q)^2 q \rceil
\label{eqn:thrm_liang_r_q_ceil}
\end{equation}
for \(q\geq Q\).
Next, we put
\[
  x = \sqrt{q}, \;\;\;  y = (1-\ve_q) \sqrt{q}
\]
and choose
\begin{equation}
     \ell = \ell(q) = \lfloor x^{1/3} \rfloor
\label{eqn:thrm_liang_ell_q_loglog_x}
\end{equation}
for \(q\geq Q\).
Finally, we define 
\begin{equation}
    k = k(q) = \left\lfloor\tfrac{1}{4}(\ell-2)^2-(\ell-2)\right\rfloor - 2
\label{eqn:thrm_liang_k_q_ell_q}
\end{equation}
for \(q\geq Q\).

We note that the three positive integers \(r,\ell,k\in\nfld\),
depending on \(q\), as chosen in the above
satisfy condition~1) and condition~2)  of Theorem~\ref{thrm:liang_applicable}
for \(q\geq Q\).
We are going to prove that condition~3) of Theorem~\ref{thrm:liang_applicable}
also holds for the chosen \(r,\ell,k\in\nfld\) for \(q\geq Q\).

The logarithmic integral function is defined as
\(\Li(x) =\int_2^x\frac{1}{\log t} dt\) for \(x\geq 2\)
\cite[p.~180]{bib:montgomery06a}.
Note that for \(2\leq y < x\), \(\Li(x)-\Li(y) \geq (x-y)/\log x\).
In \cite{bib:bennett18a}, Bennett, Martin, O’Bryant, and Rechnitzer have given
explicit bounds for primes in arithmetic progressions and specifically the following result
(see \cite[pp.~433--434]{bib:bennett18a}).

\begin{PP}
If \(x \geq 10^3\), then
\begin{equation}
\left|\pi(x;4,3) -\frac{1}{2}\Li(x)\right| < 0.53\frac{x}{(\log x)^2}.
\label{eqn:bennett_explicit}
\end{equation}
\label{prop:bennett_explicit}
\end{PP}

We are going to prove that the following cardinality of a subset
\[
   N_q
\dfnt
 \left|\{ p\in\Pfld \mid r \leq p^2 \leq q,\; p>p_\ell,\mbox{ and } p \equiv 3 \!\!\!\! \pmod {4} \}\right|
\]
satisfies \(N_q \geq 2k\) for all \(q\geq Q\).

It is computed that, for \(q \geq Q=\exp (29)\),
\(
     \ell 
= \lfloor x^{1/3} \rfloor
= \lfloor q^{1/6} \rfloor
\geq \lfloor \exp (29/6) \rfloor
= 125
\).
Then, we can check that, for all \(q \geq Q\), 
\begin{equation}
      \sqrt{r} 
   \geq (1-\ve_q) \sqrt{q}
   = (1-\ve_q) x
   > 0.6745x
   > 24x^{1/3}\log(x^{1/3})
   \geq 24\ell\log\ell
   > p_\ell,
\label{eqn:the_last_inequality}
\end{equation}
where the last inequality follows from \cite[Theorem~4.7]{bib:apostol76a}.

Therefore, we get
\[
   N_q
 = \left|\{ p\in\Pfld \mid r \leq p^2 \leq q,\; p \equiv 3 \!\!\!\! \pmod {4} \}\right|
 = \left|\{ p\in\Pfld \mid y^2 \leq p^2 \leq x^2,\; p \equiv 3 \!\!\!\! \pmod {4} \}\right|.
\]
Hence, 
\[
    N_q \geq \pi(x;4,3)-\pi(y;4,3)
\]
for all \(q \geq Q\).
According to Proposition~\ref{prop:bennett_explicit}, we obtain
\[ 
       N_q 
   \geq \pi(x;4,3)-\pi(y;4,3)
   \geq \frac{\Li(x)}{2}-\frac{\Li(y)}{2}-\frac{0.53x}{(\log x)^2}-\frac{0.53y}{(\log y)^2}
   \geq \frac{x-y}{2\log x} - \frac{1.06x}{(\log y)^2}
\]
where \(x = \sqrt{q}\) and \(y = (1-\ve_q) x\).

Now, we can verify that, for all \(q \geq Q\),
\[
       N_q 
   \geq \frac{\ve_q x}{2\log x} - \frac{1.06x}{(\log y)^2}
   \geq \frac{x}{2^{4/3}(\log x)^{4/3}} - \frac{1.06x}{(\log x+\log 0.6745)^2}
   \geq \frac{1}{2}x^{2/3}
   \geq \frac{1}{2}\ell^2
   \geq 2k.
\]

Therefore, the three positive integers \(r,\ell,k\in\nfld\)
as chosen in (\ref{eqn:thrm_liang_r_q_ceil}), (\ref{eqn:thrm_liang_ell_q_loglog_x}), and (\ref{eqn:thrm_liang_k_q_ell_q}),
depending on \(q\), satisfy all three conditions
in Theorem~\ref{thrm:liang_applicable} for all \(q \geq Q\).
Hence, it follows from (\ref{eqn:Rq_general_bound}) that
\begin{equation}
        R_q(\delta) 
   \geq (1-\delta)\frac{\log r}{\log q} - \frac{\log\sqrt{D}}{k\log q} 
   \dfnt R_\mathrm{NFC}(\delta,q)
\label{eqn:R_q_delta_apply}
\end{equation}
for all \(q \geq Q\) and \(0 < \delta < 1\), where \(D = 4 p_1\cdots p_\ell\).

Recall that the \(q\)-ary Gilbert--Varshamov bound is given by
\[
    R_\mathrm{GV}(\delta,q)
   = 1 - \delta\log_q(q-1)-\delta\log_q\frac{1}{\delta}-(1-\delta)\log_q\frac{1}{1-\delta}
\]
for all \(q \geq 2\) and \(0 < \delta < 1-q^{-1}\).
By (\ref{eqn:R_q_delta_apply}), it suffices to verify that, for all \(q \geq Q\),
\begin{equation}
      R_\mathrm{NFC}(1/2,q)
    >
      R_\mathrm{GV}(1/2,q),
\label{eqn:R_q_improve}
\end{equation}
which implies that the \(q\)-ary Gilbert--Varshamov bound can be improved for all \(q \geq Q\).

For all \(q \geq Q\), we have
\begin{equation}
      R_\mathrm{GV}(1/2,q)
  = \frac{1}{2} + \left( \frac{1}{2}\log \frac{1}{1-q^{-1}} - \log 2 \right) \frac{1}{\log q}
  < \frac{1}{2} - 0.6839\frac{1}{\log q}.
\label{eqn:GV_1over2}
\end{equation}
On the other hand, noting that
\[
         \frac{1}{2}\log\frac{q}{r}
    \leq \frac{1}{2}\log\frac{q}{(1 - \ve_q)^2 q}
       = \log\frac{1}{1 - \ve_q}
       <  \log\frac{1}{0.6745}
       < 0.3938
\]
for all \(q \geq Q\), it follows from (\ref{eqn:R_q_delta_apply}) that
\begin{equation}
      R_\mathrm{NFC}(1/2,q)
  = \frac{1}{2} - \left( \frac{1}{2}\log\frac{q}{r} + \frac{\log D}{2k} \right) \frac{1}{\log q}
  > \frac{1}{2} - \left( 0.3938 + \frac{\log D}{2k} \right) \frac{1}{\log q}
\label{eqn:NFC_1over2}
\end{equation}
for all \(q \geq Q\).

Therefore, in order to prove (\ref{eqn:R_q_improve}) for all \(q \geq Q\),
it suffices to verify that, for all \(q \geq Q\),
\begin{equation}
      \frac{\log D}{2k} 
    \leq 0.6839 - 0.3938
    =   0.2901
\label{eqn:D2k}
\end{equation}
where \(D = 4 p_1\cdots p_\ell\) and  \(k = \left\lfloor\tfrac{1}{4}(\ell-2)^2-(\ell-2)\right\rfloor - 2\).

In terms of Chebyshev's \(\vartheta\)-function \cite[p.~75]{bib:apostol76a}
defined by
\(
     \vartheta(x) = \sum_{p\leq x} \log p
\)
for \(x\geq 2\),
the verification of (\ref{eqn:D2k}) for all \(q \geq Q\) is equivalent to proving
\begin{equation}
    \vartheta(p_\ell) 
 \leq 
    -\log 4 +0.5802 \left[ \left\lfloor \tfrac{1}{4}(\ell-2)^2-(\ell-2) \right\rfloor - 2 \right]
\label{eqn:Chebyshev_theta}
\end{equation}
for all \(q \geq Q\).

Since \(q \geq Q=\exp (29)\), we have
\(\ell = \lfloor q^{1/6} \rfloor \geq 125\) and hence \(p_\ell \geq p_{125}=691\).
Thus, by applying an explicit bound for \(\vartheta\)-function given in \cite[p.~211]{bib:rosser41a},
we have
\begin{equation}
      \vartheta(p_\ell) 
   <
     \left(1+\frac{3}{\log p_\ell}\right) p_\ell.
\label{eqn:theta_function}
\end{equation}
Moreover, by using an explicit bound given in \cite[p.~433]{bib:bennett18a},
we have
\begin{equation}
      \frac{p_\ell}{\log p_\ell} < \pi(p_\ell)=\ell.
\label{eqn:pi_function}
\end{equation}

Noting the last inequality in (\ref{eqn:the_last_inequality}), i.e.,
\(p_\ell < 24\ell\log\ell \),
it follows from (\ref{eqn:pi_function}) that
\[
   p_\ell 
< \ell \log p_\ell 
< \ell \log \ell + \ell \log\log p_\ell
< \ell \log \ell + \ell \log\log \ell^2.
\]
Hence, by using (\ref{eqn:theta_function}) and (\ref{eqn:pi_function}), we get
\[
        \vartheta(p_\ell)
    <  3\ell +p_\ell
    < (3+\log 2)\ell +\ell\log\ell +\ell\log\log\ell.
\]

Therefore, in order to prove (\ref{eqn:Chebyshev_theta}) for all \(q \geq Q\),
it is sufficient to verify that
\[
    3.7\ell +\ell\log\ell +\ell\log\log\ell
 \leq
  -1.39 +0.58\left[\tfrac{1}{4}(\ell-2)^2-(\ell-2) - 3\right]
\]
for all \(\ell \geq 125\). This can be easily checked.
The proof of Theorem~\ref{thrm:liang} is complete.
\eof


\section{Acknowledgment}
\label{sec:ack}
The author would like to thank Alexander Barg for many useful conversations and
suggestions. Specifically, the author is grateful for his pointing out the reference
\cite{bib:zinoviev85a} and providing the author some guidance in understanding the result therein.



\indent
\mbox{} \\
\indent
{\small Louisiana State University, Baton Rouge, LA 70803}
\mbox{} \\
\indent
{\small \em E-mail:} xbliang@lsu.edu

\end{document}